\documentclass[a4paper]{article}

\title{Optimal control models of the goal-oriented human locomotion}

\date{\today}

\author{\textsc{Yacine
Chitour}\thanks{\scriptsize Universit\'e Paris XI, Laboratoire des
signaux et syst\`emes (L2S) Sup\'elec, 3 rue Joliot-Curie, 91190
Gif-sur-Yvette, France. \texttt{chitour@lss.supelec.fr}.}, \
\textsc{Fr\'ed\'eric Jean}\thanks{\scriptsize ENSTA, UMA, 32
boulevard Victor, 75739 Paris, France.
\texttt{Frederic.jean@ensta.fr}.}, \ \textsc{Paolo
Mason}\thanks{\scriptsize Universit\'e Paris XI, Laboratoire des
signaux et syst\`emes (L2S) Sup\'elec, 3 rue Joliot-Curie, 91190
Gif-sur-Yvette, France. \texttt{mason@lss.supelec.fr}.\newline This work was supported by the Digiteo grant {\it Congeo}} \hfill
\break }

\usepackage{graphicx,enumerate}
\usepackage[usenames]{color}
\usepackage{epsfig}
\usepackage{srcltx}

\newcommand{\al}{\alpha}
\newcommand{\vp}{\varphi}

\newcommand{\comm}{}\newcommand{\commf}{}

\newcommand{\beq}{\begin{equation}}
\newcommand{\eeq}{\end{equation}}
\newcommand{\be}{\begin{equation}}
\newcommand{\ee}{\end{equation}}
\newcommand{\bea}{\begin{eqnarray}}
\newcommand{\eea}{\end{eqnarray}}
\newcommand{\brs}{\begin{eqnarray*}}
\newcommand{\ers}{\end{eqnarray*}}
\newcommand{\ba}{\begin{array}}
\newcommand{\ea}{\end{array}}
\newcommand{\br}{\begin{eqnarray}}
\newcommand{\er}{\end{eqnarray}}

\newcommand{\Cv}{C_{\vp,\psi}}

\newcommand{\la}{\left\langle}
\newcommand{\ra}{\right\rangle}

\newcommand{\Je}{\mathcal{C}}

\def\eps{\varepsilon}

\newcommand{\lam}{\lambda}

\newcommand{\costu}{\Je(u(\cdot))}

\usepackage{amsthm,amssymb,amsbsy,amsmath,amsfonts,amssymb,amscd,mathrsfs}

\theoremstyle{plain}
\newtheorem{theorem}{Theorem}[section]

\newtheorem{lem}[theorem]{Lemma}

\newtheorem{prop}[theorem]{Proposition}
\theoremstyle{definition}

\theoremstyle{remark}
\newtheorem{rem}[theorem]{Remark}
\numberwithin{equation}{section}

\newcommand{\N}{{\mathbf N}}

\newcommand{\R}{{\mathbb{R}}}

\newcommand{\bi}{\begin{itemize}}
\newcommand{\ei}{\end{itemize}}

\newcommand{\con}{\mathcal{C}}
\newcommand{\W}{\mathcal{W}}

\begin{document}
\maketitle

\begin{abstract}
\comm{In recent papers it has been suggested that human locomotion may
  be modeled as an inverse optimal control problem. In this paradigm,
  the trajectories are assumed to be solutions of an optimal control
  problem that has to be determined. We discuss the modeling of both
  the dynamical system and the cost to be minimized, and we analyze
  the corresponding optimal synthesis. The main results describe the
  asymptotic behavior of the optimal trajectories as the target point
  goes to infinity.} 
\end{abstract}


\section{Introduction}


Ask a person walking in a empty room to leave this room through a
given door. Which path does that person choose?
The purpose of this paper is to deal with this issue, that is to
understand the goal-oriented human locomotion.

\begin{figure}[t]
\centering
\includegraphics[width=0.5\textwidth]{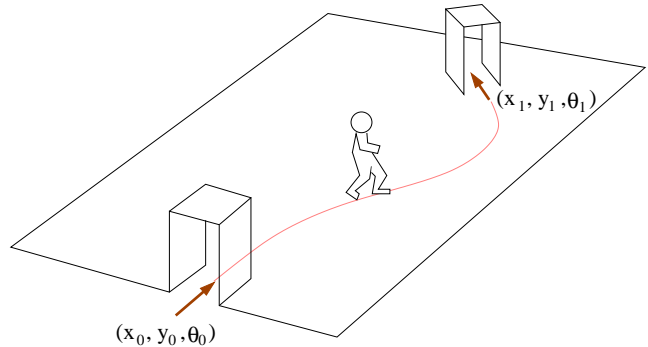}
\caption{Goal-oriented human locomotion.}
\label{f-locom}
\end{figure}

We follow the approach
initiated in~\cite{laumond2,laumond1,laumond0}. In this framework, the
locomotor trajectories lie in the simple 3-D space of both the
position $(x,y)$ and the orientation $\theta$  of the body. The
problem we want to address is then the following one. Given an initial
point $(x_0,y_0,\theta_0)$ and a final point $(x_1,y_1,\theta_1)$
(see Figure~\ref{f-locom}), which trajectory is experimentally the most
likely?

A nowadays widely accepted paradigm in neurophysiology is 
that, among all possible movements, the
accomplished ones satisfy suitable optimality criteria
(see \cite{Todorov2006} for a review). One is then led to make the assumption
that the chosen locomotor trajectory is solution of some optimal
problem, namely: minimize some integral cost
$$
\Je =  \int L(x,y,\theta, \dot x,\dot y,\dot \theta, \dots) dt
$$
among all ``admissible'' 
trajectories joining the initial point to the
final one. Two questions are in order. First, what are the
``admissible''
trajectories? It is in particular necessary to
  precise both the dynamical constraints applied to the locomotion,
  and the regularity of the trajectories (and so the order of
  derivation entering in the cost function $L$). Secondly, how to
  choose the cost function? Once the set of admissible
  trajectories is defined, the latter  question takes the form of
  an \emph{inverse optimal control problem}:  given recorded
  experimental data, infer a cost function $L$ such
  that the recorded trajectories are optimal solutions of the associated
  minimization problem.

In the theory of linear-quadratic control, the question of which
quadratic cost is minimized in order to control a linear system along
certain trajectories was already raised by R.~Kalman
\cite{Kalman1964}. Some methods allowed deducing cost functions from
optimal behaviour in system and control theory (linear matrix
inequalities, \cite{Boyd1994}) and in Markov decision processes
(inverse reinforcement learning, \cite{Ng2000}). However all these methods
have been conceived for very specific systems and they are not
suitable in the general case. A new and promising approach of the
inverse optimal control problem has been developped
in~\cite{berret,Jean} for the pointing movements of the arm. In that
approach, based on Thom transversality theory, the cost structure is
deduced from qualitative 
properties highlighted by the experimental data. Thus, starting from
the observation of inactivity intervals of the muscles during pointing
movements,  it has been proven that their presence is  a sort of necessary and
sufficient condition for the cost to be of a certain type (namely of
``absolute work'' 
type). 

The method chosen in the present paper bears some resemblance with 
~\cite{berret,Jean}. However no strong qualitative properties such as the
inactivations come out in the registered data, which prevent the
use of transversality theory. We then choose a more direct approach,
consisting of three steps: 
\begin{enumerate}

\item Modeling step: use experimental observations to define the set
  of admissible trajectories and to reduce the class of acceptable
  cost functions $L$; as a result, we obtain a class of optimal
  control problems  indexed by the cost function
  $L$. \smallskip

\item Analysis step: make a qualitative analysis of the optimal synthesis of
 the above-defined problems; the aim is mainly to exhibit properties
  characterizing the dependence on $L$ of the synthesis.\smallskip

\item Comparison step: through a numerical study  based on the
  characteristic properties exhibited above and a comparison with the
  trajectories experimentally
recorded,  determine which is the cost function $L$ which best fits the
experimental data.

\end{enumerate}

The last step is of different nature than the first two ones. It
requires the processing of a large number of registered data together
with a numerical study. In this paper we will then be concerned only
with the first two steps, the last one being the object of a
forthcoming study. \comm{Note that in the paper \cite{bayen} we already performed a similar analysis for a simplified model.} 

The modeling step will be addressed  in
Section~\ref{s-definition} and comes up with an optimal control problem {\bf (OCP)}. It is based on experimental observations together with
two modeling assumptions, \textbf{(A1)} and \textbf{(A2)}, and also on
 technical hypothesis, \textbf{(H1)}--\textbf{(H4)}.

Section~\ref{s-ocp2} is dedicated to the analysis step, all the
technical proofs being postponed to Section~\ref{s-technique}. We first show the existence of optimal solutions and 
apply the Pontryagin maximum principle to Problem {\bf (OCP)}. 
After a detailed analysis of the extremal flow, we show that the angle $\theta$ along an optimal trajectory is solution of a fourth-order differential equation $(OPT)_p$ with parameter $p$. In case the initial point and the
target are far enough one from each other, we are able to compute the asymptotic value $p^*$ of the parameter and to prove that the orbit corresponding to an optimal trajectory lies in the stable manifold of the (unique) unstable equilibrium of $(OPT)_{p^*}$.
This asymptotic behavior exactly provides the type of characteristic properties of the trajectories of human locomotion we
are looking for in order to test the adequacy of the model to registered data. The latter task will be the subject of a future work.

Moreover, based on the above asymptotic property, one can device a
numerical method for determining the initial value of the adjoint
vector and thus completely integrating the extremal flow. From the
point of view of optimal control, this geometric method is the main
technical novelty of our paper. It belongs to the realm of indirect
methods, but is clearly of 
different nature than the shooting methods.


\section{Optimal control based modeling}
\label{s-definition}

The first task is to define what is the set of admissible
trajectories. It is important here to precise the typical situation we
want to model, namely a person entering a room through one door and leaving by
another one. Hence, besides the fact that initial and final positions
and orientations are fixed, the velocity is positive at the extremities (the
person walks in and out the room). Moreover the two doors, and so
initial and final positions, are
supposed not too close one to the other.

Our modeling is based on
three important experimental
observations made in~\cite{laumond1,laumond0}.

\begin{figure}
\label{f-perp}
\centering
\includegraphics[width=0.3\textwidth]{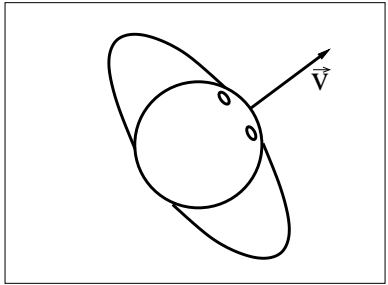}
\caption{The nonholonomic constraint.}
\end{figure}

\begin{enumerate}[$(i)$]

  \item \emph{The velocity is perpendicular to the body (see
      Figure~\ref{f-perp}).}  Namely, denoting by $\theta$ the
    orientation of the torso and by $(x,y)$ the horizontal projection
    of the center of mass, the locomotion is submitted to the
    classical nonholonomic constraint:
$$
      \dot x \sin\theta-\dot y \cos\theta=0, \qquad \hbox{i.e.} \quad \left\{\begin{array}[c]{l} \dot x= v \cos\theta \\  \dot
y=v \sin\theta \end{array}\right.
      $$
where $v=\sqrt{\dot x^2 + \dot y^2}$ is the tangential velocity.

  \item \emph{The tangential velocity has a positive lower bound: $v
      \geq a >0$.} As a consequence, all locomotion trajectories can
    be parametrized by the arc-length. We are interested here in the
    geometric curves of the human locomotion, not by their
    time-parametrized trajectories. We will then use the arc-length
    parametrization, which amounts to set $v\equiv 1$.

      It is worth to notice that, on recorded data, the tangential velocity is generally almost constant along the trajectory. Thus, up to  rescaling, arc-length and time parametrization approximately coincide.

  \item \emph{The curvature varies continuously.} That is, $\kappa
    (t)=\dot \theta (t)$ is a continuous function.
\end{enumerate}

It results from these observations that every locomotion trajectory
satisfies
$$
\dot{x}=\cos \theta, \quad
\dot{y}=\sin \theta, \quad
\dot{\theta}=\kappa,
$$
where $\kappa$ is a continuous function. Thus $\dot{x}, \dot{y}$ and
all their derivatives only depends on $\theta$ and its
derivatives. To complete the description of the set of admissible
trajectories, it only remains to precise in which functional space is
chosen the function $\theta$. Remind that our main modeling assumption
is that human locomotion is governed by optimality criteria. It is
then  mandatory to choose a functional space for the trajectories together with a class of cost functions so that one has
general results ensuring the existence of optimal solution under reasonable hypothesis. 
As for the functional space, it  must be complete and the natural choice is that of the
Sobolev spaces $W^{k,p}$. In other words, we assume that the
locomotion trajectories are solutions of the control system
$$
\begin{cases}
\dot{x}=\cos \theta, \\
\dot{y}=\sin \theta, \\
\theta^{(k)}=u,
\end{cases}
$$
where $u$ is a measurable function of finite $L^p$ norm, for a certain integer $p>0$, and $k >0$ is an integer.
As for the cost functions $L$, the highest order of derivatives of $\theta$ actually appearing in $L$ must
be equal to $k$.  The observation $(iii)$
above suggests that the order $k$ is at least equal to two. We
will make here the assumption that it is exactly equal to two.

\begin{itemize}
  \item[\textbf{(A1)}] The cost function $L =
    L(x,y,\theta,\dot{\theta},\ddot{\theta})$  explicitly depends on
    $\ddot{\theta}$, and the latter has a finite $L^p$ norm.
\end{itemize}

To sum up, the control system modeling the dynamics of the locomotion
writes as:
\begin{equation}{\label{sys1-1}}
\begin{cases}
\dot{x}=\cos \theta, \\
\dot{y}=\sin \theta, \\
\dot{\theta}=\kappa, \\
\dot{\kappa}=u,
\end{cases}
\end{equation}
where $X=(x,y,\theta,\kappa)$ belongs to $\R^2\times S^1\times \R$ and
the control $u$ is a measurable
function defined on an interval $[0,T]$, where $T>0$ depends on $u$,
and taking values in $\R$.\medskip


As for the initial and final points $X_0,X_1$,  the
spatial and angular components $(x,y,\theta)$ are given by the data of
the goal-oriented locomotion problem we consider. Two reasonable
conditions could be assumed on the curvature:
\begin{itemize}

\item either $\kappa$ is let free, the cost $\Je(u(\cdot))$ being
  minimized among
all the points $X_0,X_1$  with fixed spatial and angular components;

\item or it satisfies $\kappa = 0$ at $X_0,X_1$ (the trajectory starts and
  ends  in straight line).

\end{itemize}

This second hypothesis corresponds to  the particular experimental
setting 
considered in~\cite{laumond1}. In this
paper we will follow this assumption, however the other one
essentially leads to the same results and we will mention in
Section~\ref{s-free} how they differ.


The locomotion trajectories between two arbitrary configurations $(x_0,y_0,\theta_0)$
and $(x_1,y_1,\theta_1)$  then appear as the solutions of an optimal
control problem of the following form: minimize the cost
$$
\Je(u(\cdot))=\int_0^T L ( x,y,\theta, \kappa, u) dt
$$
among all admissible controls $u(\cdot)$ steering System~\ref{sys1-1}
from $X(0)=(x_0,y_0,\theta_0,0)$ to $X(T)=(x_1,y_1,\theta_1,0)$. Notice that the time is free in this
problem, $T$ depends on $u(\cdot)$.\medskip

The inverse optimal control problem consists in finding the proper cost function $L$. Some elementary remarks allows us to reduce the class of the candidates.

\begin{enumerate}[$(i)$]

\item The cost function is defined in an egocentric frame; the whole problem is then invariant by rototranslations, and $L$ is independent of
$(x,y,\theta)$, that is,
$$
L = L(\kappa, u).
$$

\item Turning right or left are equivalent, hence the symmetry of $L$ : $L(-\kappa, -u)= L(\kappa, u)$.

\item The more one turns, the more it costs. In other words, the
  partial functions $\kappa \mapsto L(\kappa, u)$ and $u \mapsto
  L(\kappa, u)$  are non-decreasing function of $|\kappa|$ and
  $|u|$ respectively.

\item Going in  straight line
  has a positive cost, $L(0,0)>0$, which
  is the unique minimum of $L$. Since $L$ is defined up to a
  multiplicative constant, we can impose the normalization: $L(0,0)=1$
  and $L(\kappa, u) >1$ for $(\kappa, u)\neq 0$. As a consequence, the
  cost of a straight line is its Euclidean length.

\item As already mentioned, the optimal control problem must have a
  solution, which requires some convexity properties on $L$.

\end{enumerate}

In this paper, we will restrict the class of cost functions satisfying
the following assumption:
\begin{itemize}
  \item[\textbf{(A2)}]  The cost function $L$
depends separately on  $\kappa$ and
$u$.
\end{itemize}
We will then consider a cost
\begin{equation}
\label{cost}
\Je(u(\cdot))=\int_0^T  [1+\varphi(\kappa(t))+\psi(u(t))]dt.
\end{equation}
This cost appears as a compromise between the total time $T$
(equivalently, the length to be covered) and an ``energy term''
depending separately on  $\kappa$ and $u$.

To satisfy the properties listed above we assume that the
functions $\varphi$ and $\psi$ verify the following technical hypotheses.
\bi
\item[{\bf (H1)}] $\varphi$ and $\psi$ are non negative, $\con^2$ and
  even functions defined on $\R$, and are non
  decreasing on  $\R^+$. Moreover,
  $\varphi(0)=\psi(0)=0$;
\item[{\bf (H2)}]  $\psi$ is strictly convex \comm{and $\psi''(0)>0$};
\item[{\bf (H3)}]  there exist $p>1$
and two positive constants $C,R$ such that
\begin{equation}\label{est0}
\psi(r)\geq C |r|^p, \hbox{ for every }r\in\R\hbox{ such that }|r|\geq R.
\end{equation}
\ei


\medskip

We next provide some explanations about these hypotheses.

  Hypothesis {\bf (H1)} directly results from properties
$(ii)$--$(iv)$ and is structural to the model, except for the regularity assumption. However, if the cost was not differentiable, the phenomenon of ``inactivation" would appear in the optimal trajectories (see~\cite{Jean}). Since inactivation  is not observed in the registered data, one must assume some regularity property for the cost. \comm{Note that the $\con^2$ regularity of $\vp$ and $\psi$ required in {\bf (H1)} is actually needed just for the development of the numerical method of Section~\ref{s-num-fine}, while the qualitative results of Sections~\ref{s-exis} and~\ref{s-pmp} are valid even if we assume that  $\vp$ and $\psi$ are just $\con^1$.}

\comm{The convexity and the growth condition provided by}
Hypotheses {\bf (H2)}--{\bf (H3)} are
classical assumptions which 
ensure
Property $(v)$ (see Proposition~\ref{prop-ex} and Remark~\ref{rem1} below).
\comm{The technical assumption $\psi''(0)>0$ ensures that, firstly the minimizers are locally unique; secondly  they are solutions of an ODE (of pendulum-type); and thirdly one can recover at least numerically the initial value of the adjoint vector (in an approach based on the Pontryagin maximum principle).}
It is worth to notice that {\bf (H2)}--{\bf (H3)} define open conditions
on the set of functions $\varphi$ and $\psi$ satisfying {\bf
  (H1)}. This stability property is necessary due to the physiological
nature of the cost.

The purpose of Assumption {\bf (A2)} is to 
\comm{ensure that Property $(v)$ is satisfied with} reasonable hypotheses like {\bf (H2)}--{\bf (H3)} \comm{(see, for instance, \cite{LM} for a discussion about existence of minimizers under similar assumptions).} 
      


To summarize, the optimal control problem modeling the goal-oriented human locomotion is the following one.

\begin{quote}
{\bf (OCP)} Fix an initial point $X_0=(0,0,\pi/2,0)$. For every final
point of the form $X_1=(x_1,y_1,\theta_1,0)$, for some
$(x_1,y_1,\theta_1)\in\R^2\times S^1$, find the trajectories of
\eqref{sys1-1} steering the system from $X_0$ to $X_1$ and minimizing
the cost~\eqref{cost}, where $\varphi$ and $\psi$ verify hypothesis
{\bf (H1)}--{\bf (H4)}.
\end{quote}


\section{Qualitative analysis of the optimal trajectories}
\label{s-ocp2}

\subsection{\comm{Notations}}
\label{s-notations}
In the following we will make use of the following notations.

\noindent We will always assume without loss of generality
that the difference $\al_1-\al_2$ between two angles $\al_1,\al_2\in
[0,2\pi]$ takes values in the interval $[-\pi,\pi]$. In particular
with this notation the modulus $|\al_1-\al_2|$ is a continuous
function of $\al_1,\al_2$ taking values in $[0,\pi]$.

\noindent Given a subset $S$ of $[0,T]$ we will denote its
complement in $[0,T]$ by $S^c$.

\noindent The Lebesgue measure of a set $S$ will be denoted
by $\mu(S)$.

\noindent The scalar product in $\R^2$ is denoted by
$\la\cdot,\cdot\ra$.

\noindent The symbol $B(x,r)$ indicates the ball of radius
$r$ centered at $x$.

\noindent As usual, given two subsets $A,B$ of a certain
vector space the set $A+B$ is defined as
$$A+B=\{x+y\,:\,x\in A,\,y\in B\}.$$


\subsection{Existence of optimal trajectories}
\label{s-exis}
The existence of solutions to problem {\bf (OCP)} is guaranteed by the
following result.

\begin{prop}\label{prop-ex}
For every choice of $X_0$ and $X_1$ in $\R^2\times
\mathcal{S}^1\times \R$ there exists a trajectory $\bar
X(\cdot)$ of \eqref{sys1-1}, defined on $[0,\bar T]$, associated to
some control
$\bar u(\cdot)$ and minimizing $\Je(u(\cdot))$ among all the
trajectories starting from $X_0$ and reaching $X_1$.
\end{prop}

\noindent {\bf Proof.}  Consider a minimizing sequence $u_n(\cdot):
[0,T_n]\rightarrow \R$, that is
$$\lim_{n\to\infty}\Je(u_n(\cdot))=\inf_{u(\cdot)} \Je(u(\cdot))$$
where $u_n(\cdot),u(\cdot)$ steer the system from $X_0$ to
$X_1$. Since $\Je(u_n(\cdot))$ is uniformly bounded, one easily
deduces from Hypothesis {\bf (H1)} and {\bf (H3)} that  $T_n$ and
$\|u_n\|_{L^p([0,T_n])}$ must be uniformly
bounded, and therefore, up to a subsequence, we can assume that $T_n$
converges to $\bar T$. Let us extend $u_n(\cdot)$ on $[0,\infty)$ by
  taking $u_n(t)=0$ if $t>T_n$. Then,  still up to subsequences,
  $u_n(\cdot)$ converges to some $\bar u(\cdot)\in L^p([0,\bar T])$ in
  the weak topology of $L^p([0,\bar T])$. 
Moreover, the functions $\kappa_n(t)=\int_0^t u_n(s)ds$ are uniformly
continuous by H\"older's inequality and they converge
pointwise to $\bar \kappa(t)=\int_0^t \bar u(s)ds$ since
$u_n(\cdot)\rightharpoonup \bar u(\cdot)$ in $L^p$. By
Ascoli-Arzel\`a theorem, $\kappa_n$ converge to $\bar\kappa$
uniformly on $[0,\bar T]$. Therefore, from the equation
\eqref{sys1-1}, we immediately get that the trajectories $X_n(\cdot)$
corresponding to $u_n(\cdot)$ converge uniformly to the trajectory
$\bar X(\cdot)$ corresponding to $\bar u(\cdot)$ and, since
these trajectories are uniformly equi-continuous, $\bar X(\bar T) =
\lim_{n\to\infty}  X_n (T_n)=X_1$. Moreover, we deduce that the sequence $\Big(\int_0^{T_n}\varphi(\kappa_n(t))dt\Big)_n$ converges to $\int_0^{\bar T}\varphi(\bar \kappa(t))dt$.

Let us denote $L_0(u(\cdot),T)=\int_0^T \psi(u(t))dt$.
It remains to prove that $L_0(\bar u(\cdot),\bar T)=\lim_{n\to\infty}
L_0(u_n(\cdot),T_n)$. For this purpose, we first observe that the
convexity of $\psi(\cdot)$ immediately implies the convexity in
$L^p([0,\bar T])$ of the functional $L_0(\cdot,\bar T)$. It is
well-known that the weak lower semi-continuity of a convex functional
is equivalent to its strong lower semi-continuity (see for instance \cite[Corollaire~III.8]{brezis}).
From this fact, and since $L_0(u_n(\cdot),\bar T)\leq L_0(u_n(\cdot),T_n)$, to conclude the proof of the proposition it
is enough to show that $L_0(\cdot,\bar T)$ is strongly lower semi-continuous in
$L^p([0,\bar T])$. Let us consider a sequence $(v_k)_{k\geq 1}$ in
$L^p([0,\bar T])$ converging in the $L^p$ norm to a function
$v$. \comm{In particular we can extract a subsequence $(v_{k_j})_{j\geq 1}$ converging almost everywhere to $v$ and such that $\lim_{j\to\infty} L_0(v_{k_j}(\cdot),\bar T) = \liminf_{k\to\infty} L_0(v_k(\cdot),\bar T)$.} According to Egoroff theorem there exists a sequence
$(A_h)_{h\geq 1}$ of open nested
subsets of $[0,\bar T]$ such that $\lim_{h\to\infty}\mu(A_h)=0$ and
the convergence $v_{k_j}\to v$ is uniform in $[0,\bar T]\setminus
A_h$. 
We get
\br
L_0(v(\cdot),\bar T) &=& \lim_{h\to\infty} \left(L_0(v(\cdot),\bar T)-\int_{A_h} \psi(v(\tau))d\tau\right)\nonumber\\
&=& \lim_{h\to\infty}\lim_{j\to\infty} \left(L_0(v_{k_j}(\cdot),\bar T)-\int_{A_h} \psi(v_{k_j}(\tau))d\tau \right)\nonumber\\ &\leq& \lim_{j\to\infty} L_0(v_{k_j}(\cdot),\bar T)\nonumber
\er
which concludes the proof of the proposition.
\hfill $\Box$

\begin{rem}\label{rem1}
The argument we provided is inspired from the proof of Theorem 8, p. 209 of \cite{LM}. \comm{Let us stress that
it is not difficult to show through appropriate examples that the convexity of $\psi$ required by Hypothesis~{\bf (H2)} is crucial for the existence of a minimizer.}
\end{rem}

\subsection{Application of the Pontryagin maximum principle}\label{s-pmp}
In order to apply the classical Pontryagin maximum
principle (PMP)~\cite{pon62}, one  needs to know that the optimal control
$\bar u(\cdot)$
is bounded in the $L^\infty$ topology. At the present stage of the
analysis, we do not
possess that information and we therefore must rely on more
sophisticated versions of the PMP.
For instance, one readily checks that {\bf (OCP)}
meets all the hypotheses required in Theorem~2.3 of~\cite{Vinter} and
we get the following.

\begin{prop}\label{Lp-inf}
Let $\bar X(\cdot)$ be
an optimal trajectory  for {\bf (OCP)}, defined on $[0,\bar T]$ and
associated to the
control $\bar u(\cdot)$. Then this trajectory satisfies the PMP.
\end{prop}

\comm{Let us associate to System~\eqref{sys1-1} the Hamiltonian function} 
\begin{equation}{\label{Hamil1}}
H=H(X,p,u,\nu)=p_1 \cos \theta +p_2 \sin \theta +p_3 \kappa +p_4 u
-\nu(1+\varphi(\kappa)+\psi(u)),
\end{equation}
where $p=(p_1,p_2,p_3,p_4) \in \mathbb{R}^4$ is the adjoint
vector and $\nu\in\R$.

The PMP writes as follows. Let $u(\cdot)$ be an optimal
control defined on the interval $[0,T]$ and $X(\cdot)$ the
corresponding optimal trajectory, \comm{whose existence is guaranteed by Proposition~\ref{prop-ex}}. 
Then $X(\cdot)$ is an \textit{extremal trajectory}, i.e. it satisfies
the following conditions.
There exists an absolutely continuous function
$p:[0,T] \rightarrow \mathbb{R}^4$  and \comm{$\nu\geq 0$} such that the pair $(p(\cdot),\nu)$ is
non-trivial, and such that we have:
\begin{equation}{\label{eq-Hamilton}}
\begin{cases}
\dot{X}(t)= \frac{\partial H}{\partial p}(X(t),p(t),\nu,u(t)), \\
\dot{p}(t)=-\frac{\partial H}{\partial X}(X(t),p(t),\nu,u(t)).
\end{cases}
\end{equation}
The maximization condition writes:
\begin{equation}{\label{PMP}}
H(X(t),p(t),u(t),\nu)=\max_{v \in\R} H(X(t),p(t),v,\nu) \
\mathrm{for} \ a.e. \ t\in [0,T].
\end{equation}
As the final time is free, the Hamiltonian is zero (see
\cite{Trelat1}):
\begin{equation}{\label{PMP-conseq}}
H(X(t),p(t),u(t),\nu)=0, \hspace{0.2cm} \forall t \in [0,T].
\end{equation}
The equation on the covector $p$, also called \textit{adjoint equation},
becomes:
\begin{equation}{\label{dual}}
\begin{cases}
\dot{p}_1 =0, \\
\dot{p}_2 =0,\\
\dot{p}_3 =p_1 \sin \theta - p_2 \cos \theta, \\
\dot{p}_4 =-p_3+\nu \varphi'(\kappa).
\end{cases}
\end{equation}

If $\nu\neq 0$ we can always suppose, by linearity of the adjoint equation, that $\nu=1$. In this case (resp., if $\nu=0$) a solution of the PMP is called a \textit{normal extremal} (resp., an \textit{abnormal extremal}).
It is easy to see that all optimal trajectories are normal extremals.
Indeed, if $\nu=0$, then $p_4\equiv 0$ by the maximization condition
(\ref{PMP}). From $\dot{p}_4=0$, we immediately deduce that $p_3\equiv
0$ and, from $\dot{p}_3=0$, it remains $p_1 \sin \theta - p_2 \cos
\theta\equiv 0$. From $H=0$, one also has
$p_1\cos\theta+p_2\sin\theta\equiv 0$ and thus $p_1=p_2=0$. That
contradicts the non-triviality of
$(p,\nu)$.

Consequently, Equation~(\ref{PMP-conseq}) becomes
\begin{equation}\label{ham00}
p_1\cos\theta+p_2\sin\theta+p_3\kappa+p_4u-(1+\varphi(\kappa)+\psi(u))=0.
\end{equation}

As regards the maximization condition (\ref{PMP}),  the optimal
control is given by
\begin{equation}\label{p4}
u_{opt}(t)=(\psi')^{-1}(p_4(t))\quad
\mathrm{for}  \ t\in [0,T].
\end{equation}
Note that the strict convexity and the growth condition on $\psi$
imply that $\psi'$ realizes a bijection from $\R$ to $\R$ and thus its
inverse 
is a continuous and strictly increasing function from $\R$ to $\R$.


From (\ref{dual}), we get that $p_1$ and $p_2$ are constant. \comm{Let
  us write 
 $(p_1,p_2)=\rho (\cos\phi,\sin\phi)$ for some $\phi\in [0,2\pi)$ and $\rho>0$.} Therefore
from the Hamiltonian system~\eqref{eq-Hamilton} we get that, along an
optimal trajectory, the following equation, independent of $u$,  is
satisfied for a.e. $t\in  [0,T]$ and for a suitable choice of
\comm{$(\rho,\phi)\in\R\times [0,2\pi)$},
\begin{equation}
\begin{cases}
\dot\theta =\kappa, \\
\dot\kappa=(\psi')^{-1}(p_4),\\
\dot{p}_3 = \rho \sin (\theta-\phi),\\
\dot{p}_4 =-p_3+\varphi'(\kappa)\,.
\end{cases}\label{dual1}
\end{equation}
\commf{Note that this equation implies that the function $\theta$ is
  solution of a fourth order differential equation of pendulum type
  depending on the  parameters $(\rho,\phi)$.}

Setting $(x_1,y_1)=|(x_1,y_1)|(\cos\alpha,\sin\alpha)$ for some $\alpha\in [0,2\pi)$, a careful analysis of this equation yields the following results.

\begin{prop}
There exists $R>0$ and $C>0$ such that, for every optimal trajectory  with
$|(x_1,y_1)|\geq R$ one has
$$
\| (\kappa,p_3,p_4)\|_{L^{\infty}} \leq C \qquad \hbox{and} \qquad \|
u \|_{L^{\infty}} \leq C .
$$
\label{p-bornes}
\end{prop}
\begin{rem}
The boundedness on $\| \kappa\|_{L^{\infty}}$ is valid without any hypothesis on $|(x_1,y_1)|$.
\end{rem}


\begin{prop}
For every $\eta>0$ there exists $R_{\eta}>0$ such that for every
optimal trajectory with
$|(x_1,y_1)|\geq R_{\eta}$ one has
$$
|\phi-\al|\leq\eta \qquad \hbox{and} \qquad |\rho-1| \leq \eta
$$
\label{p-asymp}
\end{prop}

\begin{theorem}
Let us associate to any extremal trajectory $X(\cdot)$ of {\bf (OCP)}
the function
$Z(t)=(\theta(t),\kappa(t),p_3(t),p_4(t))$.
Given $\nu>0$ there exist $\tau_{\nu}>0$ and
$\sigma_\nu>2\tau_\nu$ such that, for every optimal trajectory with
final time $T>\sigma_\nu$, one has
$$
|Z(t)-(\al,0,0,0)|<\nu \qquad \hbox{for }  t\in [\tau_\nu,T-\tau_\nu].
$$
\label{t-struct-lim}
\end{theorem}

The proof of these results is rather long and technical, it is
postponed to Section~\ref{s-proofs} for Propositions~\ref{p-bornes} and~\ref{p-asymp}, and Section~\ref{s-main} for Theorem~\ref{t-struct-lim}.


\subsection{Numerical study of the asymptotic behavior}
\label{s-num-fine}

The main information provided in the previous section concerning \comm{optimal trajectories, solutions of {\bf (OCP)},} can be summarized as follows: 
\bi
\item[(1)] if \comm{the norm of the spatial components $(x_1,y_1)$ of $X_1$ is large}
  enough then optimal trajectories can be decomposed in three pieces
  corresponding to time intervals $[0,\bar t]$, $[\bar t,T-\bar t]$,
  $[T-\bar t,T]$, where $\bar t$ can be thought independent of $X_1$
  and the arc of the trajectory on $[\bar t,T-\bar t]$ is
  approximately a segment (the accuracy of the approximation depends
  on the size of $\bar t$);

\item[(2)] for any optimal trajectory there exist two scalars
  $(\rho,\phi)$ such that $Z(\cdot) =
  (\theta(\cdot),\kappa(\cdot),p_3(\cdot),p_4(\cdot))$ satisfies \comm{Equation~\eqref{dual1}.}
Also, the relation
\begin{equation}
\hat H(\theta,\kappa,p_3,p_4) := \rho\cos(\theta-\phi)+p_3 \kappa+p_4 (\psi')^{-1}(p_4)-1-\vp(\kappa)-\psi\big((\psi')^{-1}(p_4)\big)=0 \label{eq-ham}
\end{equation}
holds along the trajectory.
Moreover, if $(x_1,y_1)$ is large enough, $\rho$ is close to $1$ and $\phi$ is close to the angle
$\al$ such that $(x_1,y_1)=|(x_1,y_1)|(\cos\al,\sin\al)$.
\ei

The qualitative properties stressed above do not allow
neither to understand the local behavior of optimal trajectories, in particular on the intervals $[0,\bar t],\, [T-\bar t,T]$ defined by the above Condition (1), nor to find them numerically. However they detect some non-trivial asymptotic behavior of the pair $(\rho,\phi)$ and of the initial data of \eqref{dual1}, 
for large values of $(x_1,y_1)$. The analysis carried out in this section arises from the observation that, in order to understand the asymptotic shape of the optimal trajectories on $[0,\bar t],\, [T-\bar t,T]$, it would be enough to complete the information about the initial data of \eqref{dual1}. 
Indeed, as far as the initial datum of the equation is close to its asymptotic value (if it exists) and $(\rho,\phi)$ is close to $(1,\al)$, we know, from classical continuous dependence results for the solutions of differential equations, that the solution of \eqref{dual1}. 
will in turn be close (on compact time intervals) to the solution of the \textit{asymptotic equation}
\begin{equation}
\left\{\ba{l}
\dot \theta= \kappa,\\
\dot\kappa = (\psi')^{-1}(p_4),\\
\dot p_3= \sin(\theta-\al),\\
\dot p_4=-p_3+\vp'(\kappa).
\ea
\right.\label{eq-asymp}
\end{equation}
where we take as initial value the asymptotic value of the initial data for \eqref{dual1}. 
In other words, a precise knowledge of the asymptotic behavior of such initial data, for large $(x_1,y_1)$,  would provide a tool to study numerically, through \eqref{eq-asymp}, the asymptotic shape of optimal trajectories on $[0,\bar t]$ (and, by symmetry, on $[T-\bar t,T]$).

Let us first notice that an asymptotic value for $p_4(0)$ is simply provided by evaluating \eqref{eq-ham} at time $0$, with the approximation $(\rho,\phi) = (1,\al)$. More precisely $p_4(0)$ coincides with a solution $z$ of the equation
\beq
\cos(\pi/2-\al) + z\, (\psi')^{-1}(z)-1-\psi\big((\psi')^{-1}(z)\big)=0\,.
\label{eq-vinculo}
\eeq
Since the map $\eta\mapsto \psi'(\eta)\eta - \psi(\eta)=\int_0^\eta (\psi'(\eta)-\psi'(\mu))\,d\mu$ is strictly increasing for $\eta\geq 0$, strictly decreasing for negative $\eta$ and goes to infinity for $|\eta|$ going to infinity, because of the strict convexity of $\psi$, and since $\psi\in\con^1$, we know that the previous equation has exactly one positive  solution and one negative solution. Since $u(0)=(\psi')^{-1}(p_4(0))$, this suggests the existence of two asymptotic behaviors for the  trajectories of \eqref{dual1}, 
each one corresponding to a candidate solution for {\bf (OCP)}. These two trajectories  start from $X_0$ by turning on opposite directions.

To complete the information about the asymptotic value of the initial data for \eqref{dual1} 
we need to investigate the possible values of $p_3(0)$. For this purpose we will develop below a numerical method based on the existence of a \textit{stable manifold} for \eqref{eq-asymp}.\medskip

An equilibrium for \eqref{eq-asymp} is given by
$(\theta,\kappa,p_3,p_4)=(\al,0,0,0)$ and
we know from Theorem~\ref{t-struct-lim} that, for solutions of {\bf (OCP)} with
$(x_1,y_1)$ far enough from the origin, the corresponding values of
$(\theta(\cdot),\kappa(\cdot),p_3(\cdot),p_4(\cdot))$
are close to this equilibrium on some interval $[\bar t,T-\bar t]$ for large $\bar t$ and $T$,
which suggests some stability property of the equilibrium. It is actually easy to see
that $Z_{eq}=(\al,0,0,0)$ is not a stable equilibrium of the system.
Indeed the linearized system around $Z_{eq}$ is
\beq
\dot Z = J(Z-Z_{eq})\,,\quad J = \left( \ba{cccc} 0 & 1 & 0 & 0 \\ 0 &
0 & 0 & 1/\psi''(0)\\1 & 0 & 0 & 0\\0 & \varphi''(0) & -1 & 0
\ea\right)\,,\quad Z\in\R^4,\label{linear}
\eeq
where the matrix $J$ has exactly two eigenvalues $\lam_1,\lam_2$ with
negative real part, corresponding to some eigenvectors $v_1,v_2$,
while the other two eigenvalues $\mu_1,\mu_2$, with corresponding
eigenvectors $w_1,w_2$, have positive real part. Therefore $Z_{eq}$ is
a stable equilibrium for the linearized dynamics restricted to
$Z_{eq}+V$, where $V$ is  the two dimensional real subspace of $\R^4$
spanned by $v_1,v_2$ (notice that $v_1,v_2$ can be assumed either real
or complex conjugate).

The classical stable manifold theorem (see for instance \cite{katok}) ensures the existence of a manifold $\mathcal{W}^s$ of dimension $2$, called stable manifold, which is tangent to $V$ and which contains all the trajectories converging to the equilibrium (exponentially fast).
Note that, since the continuous function $\hat H$, with $(\rho,\phi) =
(1,\al)$, is a first integral of the dynamics \eqref{eq-asymp} and
$\hat H(Z_{eq})=0$ we have $\mathcal{W}^s \subset \hat H^{-1}(0)$.

\begin{figure}
\centering
\includegraphics[width=0.45\textwidth]{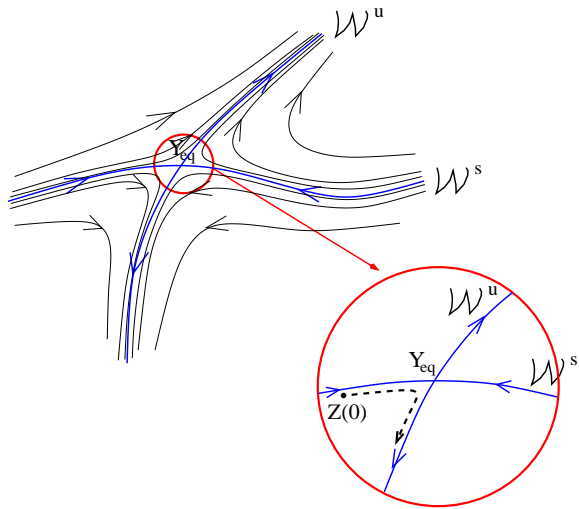}
\caption{Behavior around the stable and unstable manifolds.}
\label{ingrandim}
\end{figure}

On a small neighborhood of the equilibrium all the trajectories that
are not contained in $\mathcal{W}^s$ diverge from it exponentially
fast (see Figure~\ref{ingrandim}). Let us fix such a  neighborhood
$U$. From Theorem~\ref{t-struct-lim} we know that there exists $\bar
t$ such that, if $Z(\cdot)$ is a trajectory of \eqref{dual1} 
associated to a solution of {\bf (OCP)}, then $Z(t)\in U$ for every
$t\in [\bar t,T-\bar t]$, provided that $(x_1,y_1)$ is far enough from
the origin. 
In particular if we consider a sequence of final points $X_1^{(n)}$
for {\bf (OCP)} with spatial components $(x_1^{(n)},y_1^{(n)})=n
(\cos\al,\sin\al)$ we deduce that, for the corresponding sequence of
trajectories  $Z^{(n)}(\cdot)$, the limit $\bar Z$ of $Z^{(n)}(\bar
t)$ exists (up to a subsequence) and is contained in $\W^s$.
Continuous dependence results for the solutions of differential
equations guarantee that the limit of $Z^{(n)}(0)$ coincides with
$\bar Z(0)$, where $\bar Z(\cdot)$ is the solution of \eqref{eq-asymp}
such that $\bar Z(\bar t)=\bar Z$. In particular it must be $\bar Z(0)
= (\pi/2,0,\bar p_3,\bar p_4)$ where $\bar p_4$ satisfies
\eqref{eq-vinculo}.
\medskip

\begin{figure}
\centering
\includegraphics[width=0.4\textwidth]{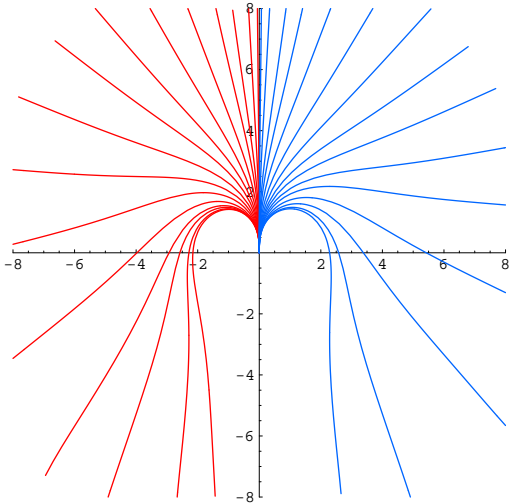}
\caption{Asymptotic behavior of optimal
trajectories with final point far from the origin.}
\label{f-asymp}
\end{figure}

The previous reasoning suggests a method to study numerically the possible values of $\bar p_3$ at time $0$.
Indeed if $U$ is small enough then $\W^s$ is well approximated by the affine space  $Z_{eq}+V$.
Consequently one can numerically look for solutions of the asymptotic equation  \eqref{eq-asymp} with
\[Z(\bar t) \in (Z_{eq}+V)\cap U\]
and such that $\theta(0)=\pi/2,\,\kappa(0)=0$. More precisely a simple
numerical method can be specified as follows. Let us fix a closed
curve $\gamma(s) = \eps\, (\cos(s) \bar v_1 + \sin(s)\bar v_2)$, where
$\bar v_1,\bar v_2$ are real vectors  spanning $V$ and $\eps$ is a
small constant (the precision of the method increases as $\eps$ goes
to zero). Since all the trajectories converging to the equilibrium
must cross this curve (in the approximation $\W^s \simeq Z_{eq}+V$) we
can recover them by following backwards in time the solutions of
\eqref{eq-asymp} starting at $Z(0) = Z_{eq}+ \gamma(s)$ up to a time
$\tilde t < 0$ such that $\kappa(\tilde t) = 0$. The candidate
approximate asymptotic trajectories we are looking for are then
determined by the values of $s$ for which, for a reasonably not too
large $|\tilde t|$ such that $\kappa(\tilde t) = 0$, we also have
$\theta(\tilde t)=\pi/2$. The value $Z(0)$ is then a candidate value
for the initial datum of a trajectory of \eqref{dual1} 
associated to a solution of {\bf (OCP)}, for large values of
$(x_1,y_1)$. Moreover this simple method allows to approximate
numerically the initial arc of such optimal trajectories (see
Figure~\ref{f-asymp} which considers the case $\vp\equiv
0,\,\psi(z)=z^2$).

An effective method to globally construct solutions of {\bf (OCP)} for large values of $(x_1,y_1)$ is the following. Define a further closed curve $\hat\gamma(\hat s) = \hat \eps\, (\cos(\hat s) \bar w_1 + \sin(\hat s)\bar w_2)$, where $ \bar w_1,\, \bar w_2$ are real vectors generating the unstable subspace $W$ (defined similarly to $V$). Assume that $\hat \eps \ll \eps\ll 1$ and consider the solutions of \eqref{dual1} 
with $\phi=0$ and starting from $Z(0)=\gamma(s)+\hat\gamma(\hat s)$, for suitable choices of $\rho,\eps, \hat \eps, s, \hat s$ such that $\hat H(Z(0))=0$.
For a fixed small enough $\eps>0$ and fixed $s\in[0,2\pi]$ it turns out that  the trajectory on intervals $[t_1,0]$, with $t_1<0$ not too large, is subjected to small variations with respect to the choice of $\hat \eps\ll \eps,\hat s\in[0,2\pi],\rho$ such that $\hat H(Z(0))=0$. In other words the trajectory approximately only depends on $\eps,s$ on the interval $[t_1,0]$.  Similarly as before, the value $s$ and the time $t_1$ can be chosen in such a way that $\kappa(t_1)=0$ and, at the same time, $\theta(t_1)= \pi/2 - \phi$, for a prescribed value $\phi$.

On the other hand for positive time the components along the stable subspace $V$ decrease exponentially as far as the components along $W$ are small so that, after a certain time, the trajectory evolves close to the unstable manifold $\W^u$ (see Figure~\ref{ingrandim}).
The dynamics at this stage essentially depends on the initial choice of  $\hat \eps$ and $\hat s$, where the first parameter determines \comm{the range of time such that the trajectory is 
confined inside $U$}, while the second one essentially determines the final angle.

\begin{figure}
\centering
\includegraphics[width=0.4\textwidth]{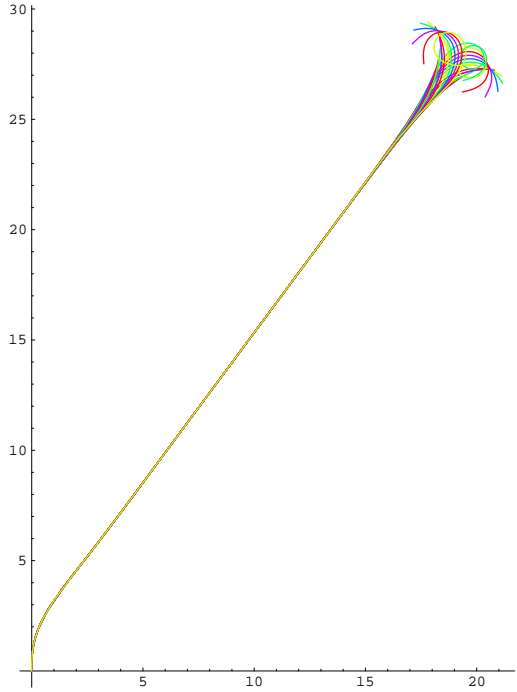}
\caption{Candidate optimal trajectories for large values of $(x_1,y_1)$.}
\label{artifici}
\end{figure}

This method gives rise, up to a rotation of an angle $\phi$ and
appropriate translations, to solutions of {\bf (OCP)}. 
Figure~\ref{artifici} depicts a set of candidate solutions of {\bf
  (OCP)} constructed by using the previous method, for a particular
choice of the angle $\phi$ and in the case $\vp\equiv 0,\psi(z) =
z^2$. 

\subsection{The case with free curvature at the extremities}
\label{s-free}

In the definition of the optimal control problem {\bf (OCP)} given in Section~\ref{s-definition}, we have chosen to take the initial and final values of the curvature equal to 0. Another reasonable condition would be to let these values free. In this case, the optimal control problem writes as follows.

\begin{quote}
$(\mathbf{\widetilde{OCP}})$ Fix an initial point $Q_0=(0,0,\pi/2)$ in coordinates $(x,y,\theta)$. \comm{For every final
 point $Q_1=(x_1,y_1,\theta_1)\in\R^2\times S^1$, find the trajectories $X(\cdot)=(x(\cdot),y(\cdot),\theta(\cdot),\kappa(\cdot))$ of
\eqref{sys1-1} such that $(x(0),y(0),\theta(0))=Q_0$ and $(x(T),y(T),\theta(T))=Q_1$, and minimizing
the cost~\eqref{cost} under hypotheses {\bf (H1)}--{\bf (H4)}.}
\end{quote}

The analysis of $(\mathbf{\widetilde{OCP}})$ leads essentially to the same results than the one of {\bf (OCP)}. The only noticeable differences are the following ones.

\begin{itemize}
  \item When applying the PMP, in addition to the Hamiltonian equations~\eqref{eq-Hamilton}, one obtain also a transversality condition on $(p(0),p(T))$, namely
      $$
      p_4(0)=p_4(T)=0.
      $$
      These conditions play the same role than the condition
      $\kappa(0)=\kappa(T)=0$ in {\bf (OCP)}.

  \item In Section~\ref{s-num-fine}, the relation~\eqref{eq-ham} on
    the Hamiltonian allows to characterize the asymptotic value
    $\kappa(0)$ (and not the one of  $p_4(0)$ through
    Equation~\eqref{eq-vinculo} as in {\bf (OCP)}), in function of
    $p_3(0)$. Indeed, since $p_4(0)=0$, one has, for the asymptotic
    values:
      $$
      \cos(\pi/2-\alpha)+p_3(0) \kappa(0)-1-\vp(\kappa(0))=0.
      $$
The numerical methods presented in Section~\ref{s-num-fine} have to be
     slightly modified in accordance with the changes above.

\end{itemize}
Note also that the proof of Lemma~\ref{lem-01} has to be modified (see
Remark~\ref{re:free}).




\section{Proofs of the main results}
\label{s-technique}

\subsection{Comparison with reference trajectories}

In order to obtain a first rough estimate of the optimal cost, we
exhibit particular trajectories  of
\eqref{sys1-1} steering the system from $X_0$ to $X_1$.

\begin{prop}
Given $\sigma>0$, a pair $(\lambda_0,\lambda_1)\in \R_+^3$ and $T\geq
2\lambda_0+2\lambda_1$ we define the control function
\beq
u(t)=\left\{\ba{rl}
-\sigma~~ & t\in[0,\lambda_0]\\
+\sigma~~ & t\in(\lambda_0,2\lambda_0]\\
0~~ & t\in(2\lambda_0,T-2\lambda_1]\\
+\sigma~~ & t\in(T-2\lambda_1,T-\lambda_1]\\
-\sigma~~ & t\in(T-\lambda_1,T]
\ea\right.\,.
\label{steps}
\eeq
Then, for every choice of
$X_1=(x_1,y_1,\theta_1,0)$ with $|(x_1,y_1)|\geq 8\sqrt{\pi/\sigma}$,
there exists a pair
$(\lambda_0,\lambda_1)\in [0,\sqrt{3\pi/\sigma}]\times
[0,\sqrt{5\pi/\sigma}]$ and $T\geq 2\lambda_0+2\lambda_1$ such that
the trajectory of
\eqref{sys1-1} with $u(\cdot)$ given by \eqref{steps} starting at
$X_0$ reaches $X_1$ at time
$T$.

\label{trajec}
\end{prop}

\begin{figure}
\centering
\input{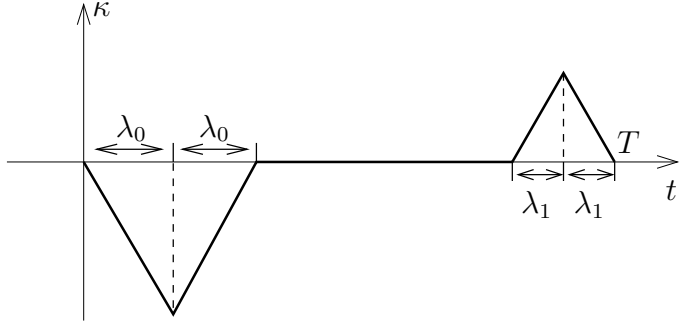}
\caption{The function $\kappa(\cdot)$ associated to the control
function of Proposition~\ref{trajec}} \label{esempio}

\end{figure}

\noindent \noindent {\bf Proof.}   First, let us observe that
$\kappa=0$ along the arc corresponding to $u=0$, which is therefore
a segment, and $\kappa(0)=\kappa(T)=0$ (see
Figure~\ref{esempio}). Also, since $\ddot \theta=u$, one can easily
check that  $\theta(t)=\theta_0-\sigma\lambda_0^2$ between
$2\lambda_0$ and $T-2\lambda_1$ and that
$\theta(T)=\theta_0-\sigma \lambda_0^2+\sigma
\lambda_1^2$. Therefore, since it must be $\theta=\theta_1$ at the
final time, and up to assuming without loss of generality that
$\theta_1-\theta_0\in [0,2\pi)$, we can assume that $\lambda_1$
is the following continuous function of $\lambda_0\in
[0,\sqrt{3\pi/\sigma}]$:
$$\lambda_1=\lambda_1(\lambda_0)=\sqrt{\frac{\theta_1-\theta_0}{\sigma}+\lambda_0^2}\,.$$
Let us denote by $\gamma_0(\cdot)$ the trajectory of \eqref{sys1-1}
corresponding to $u(\cdot)$ starting at $X_0$ at time $0$ and by
$\gamma_1(\cdot)$ the trajectory of \eqref{sys1-1} corresponding to
$u(\cdot)$ based at $X_1$ at time $T$. Our aim
is to prove that, for an appropriate choice of $(\lambda_0,T)$,
these trajectories coincide. For this purpose it is enough to prove
that their projections on the plane coincide.

\begin{figure}
 \centering
 \input{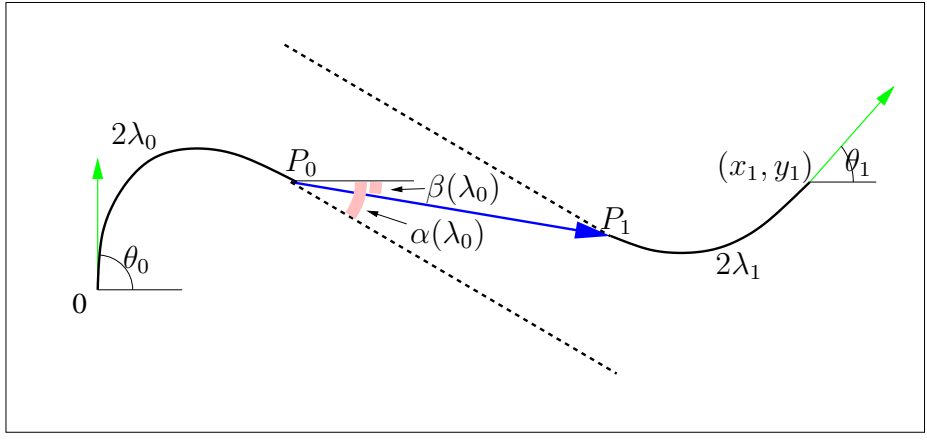}
 \caption{Proposition~\ref{trajec}, construction of the reference trajectory} \label{demo}

 \end{figure}

Let us consider the points $P_0$ and $P_1$ which are the projections
on the plane of $\gamma_0(2\lambda_0)$ and $\gamma_1(T-2\lambda_1)$
(see Figure~\ref{demo}). Since for both $\gamma_0(\cdot)$ and
$\gamma_1(\cdot)$ the angle $\theta(\cdot)$ is constantly equal to
$\al(\lambda_0):=\theta_0-\sigma\lambda_0^2$ on the interval
$[2\lambda_0,T-2\lambda_1]$, we deduce that the two curves coincide,
up to an appropriate choice of $T>0$, if and only if the angle
$\beta(\lambda_0)$ between the vector $P_1-P_0$ and the horizontal
axis is equal to $\al(\lambda_0)$ (up to a multiple of $2\pi$).

Let us observe that $\lambda_0+\lambda_1\leq
\sqrt{3\pi/\sigma}+\sqrt{5\pi/\sigma}\leq
4\sqrt{\pi/\sigma}$. Since $|P_0|< 2\lambda_0$ and
$|P_1-(x_1,y_1)|< 2\lambda_1$ when $\lambda_0\neq 0$ and
$\lambda_1\neq 0$, if $|(x_1,y_1)|\geq 8\sqrt{\pi/\sigma}$ then
$$
\begin{array}{rcl}
\la (x_1,y_1), P_1-P_0\ra & = & \la (x_1,y_1), P_1-(x_1,y_1)\ra + \la (x_1,y_1), (x_1,y_1)\ra - \la (x_1,y_1), P_0\ra \\
& \geq & |(x_1,y_1)|^2 - 2 (\lambda_0+\lambda_1) |(x_1,y_1)| \\
& \geq & 0\,.
\end{array}
$$
This implies that $\beta(\cdot)$ takes values on an interval of
length less than $\pi$. Moreover, since $P_0$ and $P_1$ are
continuous functions of $\lambda_0$, the map $\beta(\cdot)$ is also
continuous. It follows that the range of the continuous function
$\al(\cdot)-\beta(\cdot)$ is an interval whose length is larger or equal than
$2\pi$ and therefore it contains a multiple of $2\pi$, i.e. there exists $\lambda_0\in [0,\sqrt{3\pi/\sigma}]$ such that $\al(\lambda_0)=\beta(\lambda_0)$ up to a multiple of $2\pi$, as required.\hfill$\Box$

\vspace{5pt}

\noindent Comparison with the reference trajectories defined above
leads to relevant estimates as shown by the following proposition.

\begin{prop}
There exists a constant $C_{\varphi,\psi}$ only depending on
$\varphi,\psi$ such that the following holds:
if $|(x_1,y_1)|\geq
8\sqrt{\pi}$ and if $u_{opt}(\cdot)$ is an optimal control defined on
$[0,T]$ steering the system from $X_0$ to $X_1=(x_1,y_1,\theta_1,0)$
the following relations hold

\begin{equation}
|(x_1,y_1)|\, \leq \, T\, \leq\, \Je(u_{opt}(\cdot))\,\leq\,
|(x_1,y_1)|+C_{\varphi,\psi}\,. \label{bound1}
\end{equation}
Consequently,
\begin{equation}
\int_0^{T} \big(\varphi(\kappa(t))+\psi(u_{opt}(t))\big)\,dt \leq
C_{\varphi,\psi}\,.
\label{bound2}
\end{equation}

\end{prop}
\noindent {\bf Proof.} If $|(x_1,y_1)|\geq
8\sqrt{\pi/\sigma}$ we know from Proposition~\ref{trajec} that there
exists $u(\cdot)$ defined as in \eqref{steps} and steering the system
from $X_0$ to $X_1$. If $P_0,P_1$ are defined as in the proof of
Proposition~\ref{trajec} we have the following estimates
$$\begin{array}{rcl}
\Je(u(\cdot)) & = & T+2(\lambda_0+\lambda_1)\psi(\sigma)+2
\big(\int_0^{\lambda_0}\varphi(\sigma
t)dt+\int_0^{\lambda_1}\varphi(\sigma t)dt\big)\\
& \leq & T+8\sqrt{\pi/\sigma}\big(\psi(\sigma)+\|\varphi\|_{L^\infty
  ([0,\sqrt{5\pi\sigma}])}\big)\\
& = &
|P_1-P_0|+8\sqrt{\pi/\sigma}\big(1+\psi(\sigma)+\|\varphi\|_{L^\infty
  ([0,\sqrt{5\pi\sigma}])}\big)\\
& \leq &  |(x_1,y_1)|+ |P_1-(x_1,y_1)|+ |P_0|+8\sqrt{\pi/\sigma}\big(1+\psi(\sigma)+\|\varphi\|_{L^\infty ([0,\sqrt{5\pi\sigma}])}\big)\\
& < & |(x_1,y_1)|+8\sqrt{\pi/\sigma}\big(2+\psi(\sigma)+\|\varphi\|_{L^\infty ([0,\sqrt{5\pi\sigma}])}\big)\,.
\end{array}$$
Thus, by choosing $\sigma=1$ and since $\Je(u_{opt}(\cdot))\leq\Je(u(\cdot))$, we deduce the
explicit bounds \eqref{bound1}, \eqref{bound2}, with $C_{\varphi,\psi}=8\sqrt{\pi}\big(2+\psi(1)+\|\varphi\|_{L^\infty ([0,\sqrt{5\pi}])}\big)\,.$\hfill$\Box$

\begin{rem}
\label{u00}
For every $\eps>0$ and every optimal control $u_{opt}$ defined on $[0,T]$, let $\mathcal{U}_{eps}$ be the subset of $[0,T]$ given by
$$\mathcal{U}_{\eps}=\{t\in[0,T]\,:\,|u_{opt}(t)|\geq\eps\}\,.$$
From Equation~\eqref{bound2} and the strict convexity of $\psi$, we deduce that for every $\eps>0$ there exists a positive constant $C_{\vp,\psi}$ such that for every $u_{opt}$ defined on $[0,T]$, $\mu(\mathcal{U}_{\eps})\leq C_{\vp,\psi}$.
\end{rem}

As a consequence of \eqref{bound2} we easily get the uniform
equicontinuity of the $\kappa$
components of the optimal trajectories, solutions of   {\bf
  (OCP)}. This is a particular case of the following lemma, that will
also be useful in the next sections.
\begin{lem}
\label{l-equic}
For every $\Gamma>0$ and $\eps>0$ there exists
$\delta_{\eps,\Gamma}>0$ such that $|s_1-s_2|\leq\delta_{\eps,\Gamma}$
implies $|\kappa(s_1)-\kappa(s_2)|\leq \eps$ for every
$[s_1,s_2]\subset [t_1,t_2]$, whenever $t_1,t_2$ and $u(\cdot)$ are
such that $\int_{t_1}^{t_2}\psi(u(s))ds\leq \Gamma$. Moreover
$\lim_{\Gamma\to 0} \delta_{\eps,\Gamma}=+\infty$.
 \end{lem}

\noindent {\bf Proof.} By using condition {\bf (H3)} and \eqref{bound2} we get
\begin{eqnarray}
|\kappa(s_2)-\kappa(s_1)| &\leq& \int_{s_1}^{s_2} |u(\tau)| d\tau\nonumber\\
&\leq & \int_{s_1}^{s_2} \big((|u(\tau)|-R)^+ +R\big) d\tau \nonumber \\
&\leq& \Big(\int_{s_1}^{s_2} \big((|u(\tau)|-R)^+\big)^p d\tau\Big)^{1/p} |s_2-s_1|^{1/q}+R|s_2-s_1|\nonumber \\
&\leq & \Big(\frac1{C}\int_{s_1}^{s_2} \psi(u(\tau)) d\tau\Big)^{1/p} |s_2-s_1|^{1/q}+R|s_2-s_1| \nonumber\\
&\leq & \left(\frac{\Gamma}{C}\right)^{1/p}
\,|s_2-s_1|^{1/q}+R|s_2-s_1|\,, \label{eq:equicont}
\end{eqnarray}
where, for a real valued function $f(\cdot)$, we define
$f^+(t)=\max\{0,f(t)\}$ and $q$ is such that $1/p+1/q=1$. To conclude
the proof of the lemma it is enough to observe that {\bf (H3)}
actually holds for arbitrary small  $R$, provided that $C$ is also
chosen small enough.
\hfill $\Box$
\medskip

Comparisons with reference trajectories also give estimates of the
cost of pieces of trajectories that are close to a line segment.

\begin{lem}
For every $C>0$ there exists $\delta>0$ and $R>0$ large enough such
that the following holds.
Let $W_0=(\bar x_0,\bar y_0,\theta_0,\kappa_0)$, $W_1=(\bar
x_1,\bar y_1,\theta_1,\kappa_1)$, and set $(\bar x_1-\bar x_0,\bar
y_1-\bar y_0)=\Gamma (\cos\bar\theta,\sin\bar\theta)$ for some
$\Gamma>0$ and $\bar \theta\in[0,2\pi]$.
Then, if $|\theta_i-\bar\theta|<\delta$,
$|\kappa_i|<\delta$ for $i=0,1$ and $\Gamma\geq R$, any optimal
trajectory  connecting $W_0$ to $W_1$ satisfies $\costu \leq |(\bar
x_1-\bar x_0,\bar y_1-\bar y_0)|+C$. \label{innominato}
\end{lem}

\noindent {\bf Proof.}  The proof of the lemma relies on the
construction of a special trajectory  satisfying the hypotheses of the lemma, for suitable values of $\delta$ and $R$, and such that $\costu \leq
|(\bar
x_1-\bar x_0,\bar y_1-\bar y_0)|+C$. Consider
control functions of the form:
$$u(t)=\left\{\ba{rl}
\delta_0~~ & t\in[0,\tau_1]\,,\\
-\delta_0~~ & t\in(\tau_1,\tau_2]\,,\\
0~~ & t\in(\tau_2,T-\tau_3]\,,\\
\delta_1~~ & t\in(T-\tau_3,T-\tau_4]\,,\\
-\delta_1~~ & t\in(T-\tau_4,T]\,, \ea\right.$$ for some choices of
$\tau_i>0\,,\ i=1,2,3,4$, of $T>0$ and with $\delta_i
\in\{\delta,-\delta\}$ for $i=0,1$. Let us consider the two trajectories
$\gamma^0(\cdot)=(x^0(\cdot),y^0(\cdot),\theta^0(\cdot),\kappa^0(\cdot))$
and
$\gamma^1(\cdot)=(x^1(\cdot),y^1(\cdot),\theta^1(\cdot),\kappa^1(\cdot))$
corresponding to $u(\cdot)$ and such that $\gamma^0(0)=W_0$ and
$\gamma^1(T)=W_1$, respectively.
If we suppose that $\tau_2=2\tau_1+\kappa_0/\delta_0$ then a simple
computation shows that  $\kappa^0(t)=0$ on $[\tau_2,T-\tau_3]$.
We compute the value of $\theta(\cdot)$ on $[\tau_2,T-\tau_3]$
$$\theta(t)=\theta_0+\delta_0\Big[ \Big(\tau_1+\frac{\kappa_0}{\delta_0}\Big)^2-\frac{\kappa_0^2}{2\delta_0}\Big]\,,\quad t\in [\tau_2,T-\tau_3]\,,$$
and we observe that, for every fixed $\theta_0,\kappa_0$ satisfying
$|\theta_0-\bar\theta|<\delta$ and $|\kappa_0|<\delta$ and for every
value $\tilde\theta\in [\bar\theta-\delta,\bar\theta+\delta]$ there
exists $\delta_0\in \{-\delta,\delta\}$ and $\tau_1\in [0,1+\sqrt{5/2}]$
such that $\tau_2=2\tau_1+\kappa_0/\delta_0\in [0,1+\sqrt{10}]$,
$\theta(t)=\tilde\theta$. Therefore we can associate to each value
$\tilde\theta\in [\bar\theta-\delta,\bar\theta+\delta]$ a value of
$\tau_2$ and, as a  consequence, we can construct a continuous map associating $\tilde\theta\in [\bar\theta-\delta,\bar\theta+\delta]$ to a point
$(x^0(\tau_2),y^0(\tau_2))$.

Let us now consider the trajectory $\gamma^1(\cdot)$. If we set
$\tau_3=2\tau_4+\kappa_1/\delta_1$ we have that $\kappa^1(t)=0$ on
$[\tau_2,T-\tau_3]$. Moreover we have
$$\theta(t)=\theta_1-\delta_1\Big[ \Big(\tau_4+\frac{\kappa_1}{\delta_1}\Big)^2-\frac{\kappa_1^2}{2\delta_1}\Big]\,,\quad t\in [\tau_2,T-\tau_3]\,.$$
Again, it is possible to associate to each value $\tilde\theta\in
[\bar\theta-\delta,\bar\theta+\delta]$ corresponding values of
$\delta_1\in \{-\delta,\delta\}$ and $0\leq\tau_4\leq\tau_3\leq
1+\sqrt{10}$ in such a way that $(x^1(T-\tau_3),y^1(T-\tau_3))$
varies continuously with respect to $\tilde\theta$.

We want now to prove that $\gamma^0(\cdot)$ and $\gamma^1(\cdot)$
coincide, up to the choice of $T$, for a suitable value of
$\tilde\theta$. Since it is easy to see that
$(x^0(\tau_2),y^0(\tau_2))=(\bar x_0,\bar y_0)+\tau_2
(\cos\bar\theta, \sin\bar\theta)+O_1(\delta)$ and
$(x^1(T-\tau_3),y^1(T-\tau_3))=(\bar x_1,\bar y_1)-\tau_3
(\cos\bar\theta, \sin\bar\theta)+O_2(\delta)$ with $|O_i(\delta)|<M\delta$
for some universal constant $M$ we have that
\begin{eqnarray}
(x^1(T-\tau_3),y^1(T-\tau_3))-(x^0(\tau_2),y^0(\tau_2))&=&(\Gamma-\tau_2-\tau_3)(\cos\bar\theta,\sin\bar\theta)-O_1(\delta)+O_2(\delta)\nonumber\\
&=&\Gamma'(\cos\theta',\sin\theta')\nonumber
\end{eqnarray}
for some $\theta'$ with $|\theta'-\bar\theta|<\delta$ and $\Gamma'>0$, provided that $R$ is large
enough (independently of $\delta$). In particular $\theta' $ can be thought as a continuous
function of $\tilde\theta$, and we conclude that $\theta'-\tilde\theta=0$
for some value of  $\tilde\theta\in[\bar\theta-\delta,\bar\theta+\delta]$. This implies that $\gamma^0(\cdot)=\gamma^1(\cdot)$, up to choosing $T=\Gamma'+\tau_2+\tau_3$.

We have therefore constructed a trajectory
$(x(\cdot),y(\cdot),\theta(\cdot),\kappa(\cdot))=\gamma^0(\cdot)=\gamma^1(\cdot)$
corresponding to $u(\cdot)$ and connecting $W_0$ to $W_1$.

Let us estimate the cost corresponding to this trajectory. We have
that \[\int_0^T \big(\vp(\kappa(t))+\psi(u(t))\big)dt<2(1+\sqrt{10})\big(\|\vp\|_{L^\infty\left([0,\delta(2+\sqrt{5/2}\,)]\right)}+\psi(\delta)\big)\,,\] where the right-hand side can be made arbitrarily small by appropriately choosing $\delta$. It remains to compare $T$
with the difference $|(\bar x_1-\bar x_0,\bar y_1-\bar y_0)|$. We
know that $|(\bar x_1-\bar x_0,\bar y_1-\bar y_0)|=\int_0^T
\cos(\theta(t)-\bar\theta)dt$. Since $|\theta(t)-\bar\theta|\leq
M\delta$ for a suitable $M>0$ we have that $\int_0^{\tau_2}
\cos(\theta(t)-\bar\theta)dt>\tau_2(1-M^2\delta^2/2)$ and that
$\int_{T-\tau_3}^T
\cos(\theta(t)-\bar\theta)dt>\tau_3(1-M^2\delta^2/2)$. Moreover
$\int_{\tau_2}^{T-\tau_3} \cos(\theta(t)-\bar\theta)dt$ is the
projection of the segment between $P_1=(x(\tau_2),y(\tau_2))$ and
$P_2=(x(T-\tau_3),y(T-\tau_3))$ on the line $l$ connecting $(\bar
x_0,\bar y_0)$ to $(\bar x_1,\bar y_1)$. In particular since the
distance among the points $P_i$ and the line $l$ is bounded by
$M\delta$ for a suitable $M>0$, it is easy to verify that if
$T-\tau_3-\tau_2>4M^2$ then the difference among the length of the
segment between $P_1$ and $P_2$ and the length of its projection on
$l$ is bounded by $\delta^2$. In other words $\int_{\tau_2}^{T-\tau_3}
\cos(\theta(t)-\bar\theta)dt>T-\tau_3-\tau_2-\delta^2$.
By choosing $\delta$ small enough we then conclude the proof of the lemma.

$\hfill\Box$


\subsection{Some preliminary lemma}
\label{fund-lem}

Let $\alpha\in [0,2\pi)$ be such that
$(x_1,y_1)=|(x_1,y_1)|(\cos\alpha,\sin\al)$ and let us write as
 $(p_1,p_2)=\rho (\cos\phi,\sin\phi)$, for some $\phi\in [0,2\pi)$,
    the first two components of the covector associated to an optimal
    trajectory and by $\theta(\cdot)$ the corresponding angle. Note
    that the evolution of $p_3$ is described by the equation
 \begin{equation}
 \label{servira?}
 \dot p_3(t) = \rho\sin(\theta(t)-\phi)\,.
 \end{equation}

We have the following lemma.
\begin{lem}
For every $\eps>0$ there exists $\mathcal{T}_\eps>0$ such that,
for every optimal trajectory, one has  $\mu(J_\eps)\leq
\mathcal{T}_\eps$, where
the set $J_\eps$ is defined as
$$J_\eps=\{\tau\in [0,T]\,:\,|\al-\theta(\tau)|\geq\eps\}\,.$$
\label{lem2}
\end{lem}

\noindent {\bf Proof.}  Let us first note that $|(x_1,y_1)|=\la
(\cos\al,\sin\al),(x_1,y_1)\ra$. Since $(x_1,y_1)=\int_0^{T}
(\cos\theta(t),\sin\theta(t))dt$ we therefore get
$$ |(x_1,y_1)|=\int_0^{T} \cos (\al-\theta(t))dt\,. $$
In particular if $\eps$ is small enough we have
$$|(x_1,y_1)|=\int_{J_\eps} \cos (\al-\theta(t))dt  +\int_{J_\eps^c}
\cos (\al-\theta(t))dt < \Big(1-\frac{\eps^2}4\Big)\mu(J_\eps)+\mu(J_\eps^c)=
T-\frac{\eps^2}4 \mu(J_\eps).$$
Since from the uniform bound \eqref{bound1} we have
$T-C_{\vp,\psi}\leq |(x_1,y_1)|$ for $|(x_1,y_1)|$ large enough the
lemma is proved with $\mathcal{T}_\eps=4C_{\vp,\psi}/\eps^2$.\hfill$\Box$

\begin{lem}
Fix $\eta>1$ and $\bar k>0$ and, for any optimal trajectory and any
$k\geq \bar k$ define
$S_k=\{t\in [0,T]\,:\,|\kappa(t)|\in [k,k\eta]\}$. Then there exists
$C_{\vp,\psi}>0$, independent of $k$, such that $\mu(S_k)\leq
C_{\vp,\psi}$.
\label{lem-00}
\end{lem}

\noindent {\bf Proof.}
Note that $S_k$ is the union of two sets $S^+_k$ and $S^-_k$ defined respectively by the sign of $\kappa$.

For the sequel, we only provide estimates for $S^+_k$ since the same ones hold true for $S^-_k$.
For $\delta>0$, let $K_{\delta,k}$ be the union of all disjoint subintervals of $S^+_k$  of length at least $\delta$. We first
prove that there exists $C_{\vp,\psi}^{\delta}>0$ such that $\mu(K_{\delta,k})\leq C_{\vp,\psi}^{\delta}$ for any optimal
trajectory and $k\geq \bar k$.

Assume without loss of generality that $\delta<2\pi/\bar k$ and fix $\bar\eps=\delta\bar k/3$. Notice that for every
measurable set $S\subset S^+_k$ such that $\int_{S} \dot\theta(t)dt=\int_{S} \kappa(t)dt=l$ for some $l>0$ one has
$\mu(S)\geq l/(k\eta)$.

It turns out that if $[t_1,t_2]\subset S^+_k$ is such that
$\int_{t_1}^{t_2} \kappa(t)dt=3\bar\eps$ then $\mu([t_1,t_2]\cap J_{\bar\eps})\geq \bar\eps/(k\eta)$, where $J_{\bar\eps}$
is defined as in Lemma~\ref{lem2}. Clearly $t_2-t_1\leq 3\bar\eps/k\leq 3\bar\eps/\bar k=\delta$. Let now $I\subset S^+_k$ be
an interval such that
$\mu(I)\geq \delta$ and $m\geq 1$ be the largest integer such that $3\bar\eps m/k\leq \mu(I)$. Since
$\mu(I)\leq 3\bar\eps (m+1)/k\leq 6\bar\eps m/k$ it turns out that
$$\mu(I\cap J_{\bar \eps})\geq \frac{\bar\eps m}{k\eta}\geq \frac{\mu(I)}{6\eta}.$$
Thus, by definition of $K_{\delta,k}$ and the previous inequality, one gets
\begin{equation}
\mathcal{T}_{\bar \eps}\geq\mu(J_{\bar \eps})\geq\mu(K_{\delta,k}\cap
J_{\bar \eps})\geq \frac{\mu(K_{\delta,k})}{6\eta}\,,\nonumber
\eeq
hence proving that the measure of $K_{\delta,k}$ is bounded by $C_{\vp,\psi}^{\delta}=6\eta\mathcal{T}_{\eps}$.

To conclude the proof of the lemma let us define $\hat S^+_k=\{t\in [0,T]\,:\,\kappa(t)\in [k/\nu,k\eta\nu]\}$ for some
$\nu>1$ and define $\hat K_{\delta,k}$ similarly as before. Let us observe that
$[k/\nu,k\eta\nu]\supset [k,k\eta]+B(0,\eps_\star)$ for some $\eps_\star>0$ and for every $k\geq \bar k$. According to
Lemma~\ref{l-equic}, there exists $\delta_\star$ such that $|\kappa(s)-\kappa(t)|\leq \eps_\star$ if $|s-t|\leq \delta_\star$
for $s,t\in [0,T]$. One deduces that  $S^+_k\subset \hat K_{\delta_\star,k}$.
Since, as proved above, $\hat K_{\delta_\star,k}$ has uniformly bounded measure, we get the conclusion.

\hfill$\Box$

\begin{lem}
There exists a constant $\Gamma_{\vp,\psi}>0$ such that for any optimal
trajectory $\|\kappa(t)\|_{L^\infty([0,T])}\leq \Gamma_{\vp,\psi}$.
Moreover, for every $\eps>0$, there exists $\mathcal{T}_\eps>0$ such that,
for every optimal trajectory, one has  $\mu(K_\eps)\leq
\mathcal{T}_\eps$, where the set $K_\eps$ is defined as
$$
K_\eps=\{\tau\in [0,T]\,:\,|u(\tau)|\geq\eps\,\mbox{ or }\
|\kappa(\tau)|\geq\eps\}\,.
$$
\label{lem-01}
\end{lem}

\noindent {\bf Proof.}  Fix $\eta>1$. Given an optimal trajectory,
define $\hat k=\|\kappa(t)\|_{L^{\infty}([0,T])}/\eta$.
Since $\kappa$ is continuous and  $\kappa(0)=\kappa(T)=0$, the inverse
image  of  $[\hat k,\hat k\eta]$ under $\kappa$ contains a segment
$[s,t]$ such that $|\kappa(s)-\kappa(t)|=\hat k(\eta-1)$. From
Lemma~\ref{lem-00} there exists a positive constant $\hat
C_{\vp,\psi}$ independent of $\hat k$ such that $|s-t|\leq \hat
C_{\vp,\psi}$. Therefore, the uniform equicontinuity characterized in
Lemma~\ref{l-equic} (\comm{or, more directly, Equation~\eqref{eq:equicont}}) yields a
uniform bound on $|\kappa(s)-\kappa(t)|$
and thus on $\hat k$. The first part of the lemma is proved.

As regard the second part of the lemma, by using Remark~\ref{u00}, it
is enough to prove that the set
$C_{\eps}:=\{\tau\in [0,T]\,:\, |\kappa(\tau)|\geq\eps\}$ has
uniformly bounded measure.
This easily follows from Lemma~\ref{lem-00} by noticing that
$C_{\eps}=\{\tau\in [0,T]\,:\, |\kappa(\tau)|\in [\eps,\eta\eps]\}$
with $\eta=\Gamma_{\vp,\psi}/\eps$.

\hfill $\Box$

\begin{rem}
\label{re:free}
In the case where the curvature at the extremities is free, the bound
on $\|\kappa(t)\|_{L^\infty([0,T])}$ only holds when $T$ is greater
than a positive constant $T_{\vp,\psi}>0$. Indeed, at the beginning
of the proof of the lemma, the existence of a segment
$[s,t]$ such that $|\kappa(s)-\kappa(t)|=\hat k(\eta-1)$ can not be
ensured for any trajectory when  $\kappa(0)$ and  $\kappa(T)$ are not
equal to 0. However such a segment exists as soon as $T$ is
greater than the constant $C_{\vp,\psi}$ given by
Lemma~\ref{lem-00}. The rest of the proof is unchanged.
\end{rem}

We need now a simple technical lemma.

\begin{lem}
For every $m\in\N$ there exist $\mu_m>0$ and $\lambda_m\in (0,1)$
such that, for every $t_1,t_2\in\R$, $C>0$ and  $f\in \mathcal{C}^m
([t_1,t_2])$ satisfying $|f^{(m)}(t)|>C$ there exists a subinterval
of $[t_1,t_2]$ of length $\lambda_m (t_2-t_1)$ where $|f(t)|>\mu_m C
(t_2-t_1)^m$. In particular we can take $\lambda_m=4^{-m}$ and
$\mu_m=4^{-\frac{m(m+1)}{2}}$. \label{lem1}
\end{lem}

\noindent {\bf Proof.}  The proof goes by induction. Clearly the
lemma is true for $m=0$. Assume that the lemma is true up to $m=\hat
m$ and let us prove it for $m=\hat m+1$.  By the inductive
hypothesis we know that $|f^{(\hat m+1)}(t)|>C$ implies that
$|f'(t)|>\mu_{\hat m} C (t_2-t_1)^{\hat m}$ on a certain subinterval
of $[t_1,t_2]$ of length $\lambda_{\hat m} (t_2-t_1)$. For
simplicity, and without loss of generality, let us identify this
subinterval with $[0,\tau]$, where $\tau= \lambda_{\hat m}
(t_2-t_1)$, and let us assume that $f$ is increasing on it.

We claim that  $|f(t)|>\frac{\mu_{\hat m} \lambda_{\hat m} C}4
(t_2-t_1)^{\hat m+1}$ either in $[0,\tau/4]$ or in $[3\tau/4,\tau]$.
Indeed if this is not true then, since $f$ is increasing, we would
have $f(\tau/4)\geq -\frac{\mu_{\hat m} \lambda_{\hat m} C}4
(t_2-t_1)^{\hat m+1}$ and $f(3\tau/4)\leq\frac{\mu_{\hat m}
\lambda_{\hat m} C}4 (t_2-t_1)^{\hat m+1}$, which implies
$f(3\tau/4)-f(\tau/4)\leq\frac{\mu_{\hat m} \lambda_{\hat m} C}2
(t_2-t_1)^{\hat m+1}= \mu_{\hat m}  C (t_2-t_1)^{\hat m}
\frac{\tau}2$ and which contradicts the hypothesis
$|f'(t)|>\mu_{\hat m} C (t_2-t_1)^{\hat m}$. The lemma is therefore
proved by taking $\lambda_{\hat m+1}=\lambda_{\hat m}/4$ and
$\mu_{\hat m+1}=\mu_{\hat m}\lambda_{\hat m}/4$. 

\hfill$\Box$

As a consequence, we obtain a bound on the $L^\infty$ norm of $(p_3,p_4)$.

\begin{lem}\label{tech0}
There exist two positive constants $C_{\vp,\psi}$ and $T_{\vp,\psi}$
such that, for every optimal trajectory defined on $[0,T]$ with
  $T>T_{\vp,\psi}$, one has $\|(p_3,p_4)\|_{L^\infty([0,T])}\leq C_{\vp,\psi}(1+\rho)$,
where we recall that $\rho:=|(p_1,p_2)|$.
\end{lem}
\noindent {\bf Proof.} Define $\gamma_l:=p_l/(1+\rho)$ for
$l=3,4$. Then one has
$$
\dot \gamma_3=\sin(\theta-\phi),\quad
\dot\gamma_4=-\gamma_3+\vp'(\kappa)/(1+\rho),
$$
along any optimal trajectory defined on $[0,T]$. Notice that
$|\dot \gamma_3(t)|\leq 1$ and, from Lemma~\ref{lem-01},
we have $|\vp'(\kappa)|/(1+\rho)\leq K_{\vp,\psi}$ for some constant
$K_{\vp,\psi}$. From the dynamics
of $\gamma_4$, if $t_0$ is such that
$|\gamma_3(t_0)|=\|\gamma_3\|_{L^\infty([0,T])}$,
 one gets that $|\dot \gamma_4(t)|\geq
 \|\gamma_3\|_{L^\infty([0,T])}-|t-t_0|-K_{\vp,\psi}$ for $t\in[0,T]$.
From Lemma~\ref{lem1}, we deduce that, if
$\|\gamma_3\|_{L^\infty([0,T])}$ is not uniformly bounded
for $T$  large enough then equation~\eqref{bound2} cannot hold true.
Therefore we get a uniform bound on $|\dot \gamma_4|$ which, still
according to \eqref{bound2} and a
standard argument by contradiction, implies the uniform boundedness of
$\gamma_4$.  The proof of the  lemma is complete.

\hfill $\Box$

We use the previous lemma to obtain the following one.
\begin{lem}
For every $\nu>0$ there exists  $T_{\vp,\psi}>0$ such that, for any optimal trajectory defined on $[0,T]$ with
$T\geq T_{\vp,\psi}$ one has $\rho\geq 1-\nu$.\label{lem-rho}
\end{lem}

\noindent {\bf Proof.} We argue by contradiction and assume that $|(p_1,p_2)|\leq 1-\nu$.
For $\eps>0$ and $t_0\in [0,T]\setminus K_\eps$, where $K_\eps$ is defined in Lemma~\ref{lem-01},
the Hamiltonian at time $t_0$ writes
\begin{eqnarray}
H & = & p_1 \cos \theta(t_0) +p_2 \sin \theta(t_0) +p_3(t_0) \kappa(t_0) +p_4(t_0) u(t_0)
-(1+\varphi(\kappa(t_0))+\psi(u(t_0)))\nonumber\\
& \leq & -1+|(p_1,p_2)|+2C_{\vp,\psi} \eps\nonumber\\
&\leq& -\nu+2C_{\vp,\psi} \eps.\nonumber
\end{eqnarray}
Take now $\eps=\nu/4C_{\vp,\psi}$, and $T>\mu(K_{\eps})$ in the contradiction assumption, so
that $[0,T]\setminus K_\eps \neq \emptyset$. One therefore gets that
$H<0$, which contradicts the PMP. Hence the conclusion.

\hfill $\Box$

The next proposition will be crucial in order to prove our main asymptotic
results.

\begin{prop}
For every $\eta>0$ there exists $R_{\eta}>0$ such that for every
optimal trajectory with
$|(x_1,y_1)|\geq R_{\eta}$ one has $|\phi-\al|\leq\eta$ and
$\rho\leq 1+\eta$.
\label{p-lem3}
\end{prop}

\noindent {\bf Proof.}
We first establish the fact that  $|\phi-\al|\leq\eta$ if
$|(x_1,y_1)|$ is large
enough.
Assume that $|\phi-\al|>\eta$. Then, our aim is to prove that there
exists $T_{\eta}$ such that $T<T_{\eta}$. One has, for every $t,s\in[0,T]$,
$$
|\theta(t)-\theta(s)|\leq \Cv |t-s|,
$$
where $\Cv$ is provided by Lemma~\ref{lem-01}.
Set $\eps=\eta/3$ and $\delta_\eps=\eps/\Cv$. Let
$J_{\eps}$ be the set defined in Lemma~\ref{lem2} and consider the set
$$
\tilde{J}_{\eps}:=\big(J_{2\eps}+B(0,\delta_{\eps})\big)\cap
[0,T]\,.
$$
Clearly, from the definition of $\delta_{\eps}$, we have that
$\tilde{J}_{\eps}\subset J_{\eps}$. Therefore from Lemma~\ref{lem2}
we have that $\mu(\tilde{J}_{\eps})<\mathcal{T}_{\eps}$.  Moreover for
every $t\in \tilde{J}_{\eps}\setminus ([0,\delta_{\eps}]\cup
[T-\delta_{\eps},T])$ there exists a neighborhood in $[0,T]$ of
diameter at least $2\delta_{\eps}$ and completely contained inside
$\tilde{J}_{\eps}$. We deduce that $\tilde{J}_{\eps}$ is the union
of a finite number of intervals and that the number of these
intervals is bounded by $\mathcal{T}_{\eps}/(2\delta_{\eps})+1$.
As a consequence the complement $\tilde{J}_{\eps}^c$ is also the
union of at most
$n_{\eps}:=\mathcal{T}_{\eps}/(2\delta_{\eps})+2\leq\mathcal{T}_{\eps}/\delta_{\eps}$
intervals. Since $\mu(\tilde{J}_{\eps}^c)\geq T-\mathcal{T}_{\eps}\geq
\frac{T}2$ if $T$ is large enough, we get the existence of a
subinterval $I$ of $\tilde{J}_{\eps}^c$ such that
$$
\mu(I)\geq\frac{\mu(\tilde{J}_{\eps}^c)}{n_{\eps}}\geq\frac{T\delta_{\eps}}{2\mathcal{T}_{\eps}}
=:M_{\eps}T\,.
$$

In particular, since $\tilde{J}_{\eps}^c\subset J_{2\eps}^c$ implies
that $|\phi-\theta(t)|\geq |\phi-\al|-|\al-\theta(t)|>\eps$ on $I$,
we have that, if $|(x_1,y_1)|$ is far enough from the origin, one gets
according to \eqref{servira?} and Lemma~\ref{lem-rho},
\beq\label{int0}
|\dot{p}_3|>\frac12\sin\eps, \hbox{ on }I.
\eeq
By Lemma~\ref{lem1}, there exists a subinterval $I'$ of $I$ with $\mu(I')\geq M_{\eps}T/4$
and $|p_3|>N_\eps T$ on $I'$, where $N_\eps$ only depends on $\eps$. We immediately deduce
from Lemma~\ref{lem-01} \comm{and the equation of $p_4$} that $|\dot p_4|>N_\eps T/2$ on
$I'$ \comm{for $T$ large enough}. 
Again, by Lemma~\ref{lem1}, there exists a subinterval $I''$ of $I'$ with $\mu(I'')\geq M_{\eps}T/16$
and $|p_4|>N_\eps M_\eps T^2/16$ on $I''$. That contradicts Remark~\ref{u00} and we deduce that
$|\phi-\al|\leq\eta$ for $|(x_1,y_1)|$ large enough.

As for the second one, this is an immediate consequence of the above results.
Indeed, for every $\eps>0$, there exists $\bar T$ large enough such that every optimal trajectory defined on $[0,T]$ with $T\geq \bar T$ satisfies
$|\theta(t_0)-\phi|<\eps$, $|\kappa(t_0)|<\eps$ and $|u(t_0)|<\eps$, for some $t_0\in [0,T]$.
Then, from $H=0$, we get
\begin{eqnarray}
0&=&\rho\cos(\theta(t_0)-\phi)+p_3(t_0)\kappa(t_0)+p_4(t_0)u(t_0)-(1+\vp(\kappa(t_0))+\psi(u(t_0)))\nonumber\\
&\geq& \rho(1-\eps^2/2)-2 C_{\vp,\psi}(1+\rho) \eps-1-\|\vp(\cdot)\|_{L^\infty([0,\eps])}-\psi(\eps)\nonumber\\
&=& \rho(1-2 C_{\vp,\psi}\eps-\eps^2/2)-1-2 C_{\vp,\psi} \eps-\|\vp(\cdot)\|_{L^\infty([0,\eps])}-\psi(\eps),\nonumber
\end{eqnarray}
and, by taking $\eps$ small enough, the proof of the lemma is concluded.

\hfill $\Box$

\subsection{About Propositions~\ref{p-bornes} and~\ref{p-asymp}}
\label{s-proofs}

One immediately deduces from  Lemma~\ref{tech0} and
Proposition~\ref{p-lem3} that $p_3$, $p_4$ and $u$ are uniformly
bounded for $T$ large enough over all optimal trajectories. This,
together with Lemma~\ref{lem-01}, establish
Propositions~\ref{p-bornes}. Proposition~\ref{p-asymp} results from
Lemma~\ref{lem-rho} and
Proposition~\ref{p-lem3}.

At first sight it seems reasonable to conjecture that the previous results can be improved in the  following directions:
\begin{itemize}
\item[(a)] extending the uniformity results to all optimal
  trajectories, i.e. independently of the final time $T$;
\item[(b)] as the terminal point $(x_1,y_1)$ goes to infinity,  the
  corresponding optimal control $u(\cdot)$ tends to $0$.
\end{itemize}
However, the next two remarks show that it is not the case.

\begin{rem}
Disproving Conjecture (a) amounts to show that the control function $u(\cdot)$ associated to optimal trajectories
reaching points in a neighborhood of the origin is not uniformly
bounded. More precisely, we will exhibit a sequence of points
$X^{(n)}_1=\big(x^{(n)}_1,y^{(n)}_1,\theta^{(n)}_1,0\big)$
with $\lim_{n\to\infty} |(x^{(n)}_1,y^{(n)}_1)|=0$ such that the
optimal controls $u^{(n)}(\cdot)$ steering the system from $X_0$ to
$X_1^{(n)}$ satisfy $\lim_{n\to\infty}
\|u^{(n)}(\cdot)\|_{\infty}=\infty$.
For $\eps>0$, let us consider
$$u_\eps (t)=\left\{\ba{rl} \bar u_\eps~~ & t\in[0,\eps/2]\,, \\ -\bar u_\eps~~ & t\in (\eps/2,\eps]\,,\ea\right.$$
where $\bar u_\eps>0$ verifies $\eps\psi(\bar u_\eps)=m$, where $m$ is
a positive constant to be fixed later. Note that $\bar u_\eps$ tends
to infinity as $\eps$ goes to $0$ and $\bar u_\eps\leq C_0
(m/\eps)^{1/p}$ for $C_0$ independent of $\eps$ and $m$, provided that
$\eps/m$ is not too large. It implies that
$\vp_\eps:=\int_0^\eps\vp(\bar u_\eps s)ds$ tends to $0$ as $\eps$
goes to $0$.

We remark that, for every $t\in [0,\eps]$ the angle $\theta_\eps(t)$ corresponding to $u_\eps(t)$ maximizes $\theta(t)$ under the constraints $\|u\|_{\infty}\leq \bar u_{\eps}$, $\theta(0)=\pi/2$,
$\kappa(0)=0$ and $\kappa(\eps)=0$. Indeed, for every such $u(\cdot)$,
one has $\kappa(t)=\int_0^t u(s)ds \leq \int_0^t u_\eps(s)ds$, for
every $t\in[0,\eps]$.
Notice also that $ \theta_\eps(\eps)=\pi/2+\bar u_\eps
\eps^2/4$. Assuming that we start from the origin we also have that
$$y(t)=\int_0^t \sin\theta(t)dt>\int_0^t \sin\theta_\eps(t)dt=
y_\eps(t)~~\mbox{ and
}~~x(t)=\int_0^t \cos\theta(t)dt>\int_0^t
\cos\theta_\eps(t)dt=x_\eps(t)$$ where $\theta(\cdot)$ is associated
with a control $u(\cdot)$
satisfying the previous constraints and not almost
everywhere 
equal to
$u_\eps(\cdot)$. This implies that, in order to reach the point
$(x_\eps(\eps),y_\eps(\eps))$ at time $T$ with a control $u(\cdot)\neq
u_\eps(\cdot)$ satisfying the previous constraints,  it must be
$T>\eps$ and  $y(\eps)>y_\eps(\eps)$ so that there exists  $\bar
t\in[\eps,T]$ such that $\dot y(\bar t)<0$, i.e. $\theta(\bar t)\in
(\pi,2\pi)$. Let us prove that the corresponding trajectory cannot
be optimal. First we note that the total cost corresponding to $
u_\eps(\cdot)$ is $m+\eps+2\vp_\eps$. Therefore, if we assume by contradiction that
$u(\cdot)$ is optimal we must have $T<m+\eps+2\vp_\eps$ and  $\int_{0}^t \psi(u(\tau))d\tau<m+\eps+2\vp_\eps$ for every $t\in [0,T]$. From H\"older inequality, we deduce that $|\kappa(t)|<C_1 [m+\eps+2\vp_\eps]^{1/p} t^{1-1/p}$, which implies that
$|\theta(\bar t)-\pi/2|< C_1 [m+\eps+2\vp_\eps]^{2}$, for some $C_1>0$ independent of  $m$ and
 $\eps$ small enough. By taking $m=\sqrt{\pi/(3C_1)}$ we reach  a
contradiction. Therefore, any optimal control connecting $X_0$ to
$X_1=(x_\eps(\eps),y_\eps(\eps),\theta_\eps(\eps),0)$ must satisfy
$\|u\|_{\infty}>\bar u_{\eps}$ for $\eps$ small enough.
\end{rem}

\begin{rem}
We next show that Conjecture (b) is false by disproving that,
 for every $\eps>0$, there exists $R_{\eps}$ such
that  $\|u\|_{\infty}\leq \eps$ for every optimal triple
$(X(\cdot),u(\cdot),T)$ with $|(x_1,y_1)|>R_{\eps}$. Indeed, if $\al$ is
defined as in Section~\ref{fund-lem} and $|(x_1,y_1)|$ is large
enough the previous results say that $p_2$ is arbitrarily close to
$\sin\al$, and thus different from $1$ in general.  On the other hand,
$H\equiv 0$ and $\theta(0)=\pi/2$ imply $u(0)\psi'(u(0))-\psi(u(0))=1-p_2$ and therefore $u(0)$ is not close to $0$ in general.
\end{rem}

\subsection{Proof of Theorem~\ref{t-struct-lim}}
\label{s-main}

We will prove the theorem by showing separately that the functions
$\theta(\cdot)-\al,\kappa(\cdot),p_3(\cdot),p_4(\cdot)$ can be made
arbitrarily small by choosing large enough $\tau_\nu,\sigma_\nu$.

From Lemmas~\ref{lem2} and~\ref{lem-01} we have that, given $\delta'>0$
and if $T>2\mathcal{T}_{\delta'}$ for a suitably large $\mathcal{T}_{\delta'}$,  there exists
$t_0\in[0,\mathcal{T}_{\delta'}]$ and  $t_1\in[T-\mathcal{T}_{\delta'},T]$
such that, $|\al-\theta(t_i)|<{\delta'}$ and
$|\dot\theta(t_i)|<{\delta'}$. If we set $W_0=X(t_0)$ and $W_1=X(t_1)$
and we let $\bar\theta$ be as in Lemma~\ref{innominato} then, if $T$
is large enough (depending only on ${\delta'}$), we have that
$|\al-\bar\theta|<{\delta'}$ and the hypotheses of
Lemma~\ref{innominato} are satisfied with $\delta=2{\delta'}$.  
From Lemma~\ref{innominato}, and for a fixed $C>0$, we have that the optimal control
$u(\cdot)$ must be such that $\int_{t_0}^{t_1}\psi(u(t))dt\leq C$, if $\delta$ is small enough.

Let us fix $\eps>0$. If $C$ is small enough then, from Lemma~\ref{l-equic}, we can assume that
$|s_1-s_2|\leq 4\mathcal{T}_\eps$ implies
$|\kappa(s_1)-\kappa(s_2)|\leq
\tilde\eps:=\frac{4\eps}{\mathcal{T}_\eps}$, where $\mathcal{T}_\eps$
is defined by Lemma~\ref{lem2}. By contradiction it is then easy to
see that $\|\kappa(\cdot)\|_{L^\infty([t_1,t_2])}\leq2\tilde\eps$ if
$T$ is large enough. Indeed, otherwise we would have
$|\kappa(t)|>\bar\eps$ on an interval of length larger than
$4\mathcal{T}_\eps$ and then, as a consequence of Lemma~\ref{lem1}, we
would get $|\theta(t)-\bar \theta|>\eps$ on an interval of length
larger than $\mathcal{T}_\eps$ contradicting Lemma~\ref{lem2}.
Moreover, since
$\|\kappa(\cdot)\|_{L^\infty([t_1,t_2])}\leq2\tilde\eps$ we have
$\|\theta(\cdot)-\bar \theta\|_{L^\infty([t_1,t_2])}\leq
\eps+2\tilde\eps/\mathcal{T}_\eps$, again as a consequence of
Lemma~\ref{lem2}.

It remains to prove that
$\|(p_3(\cdot),p_4(\cdot))\|_{L^\infty([t_1,t_2])}$ can be made
arbitrarily small by an appropriate choice  of
$\tau_\nu,\sigma_\nu$. For this purpose it is enough to observe that
$\|\dot p_3\|_{L^\infty([t_1,t_2])}$ can be made arbitrarily small,
and, again, a simple argument by contradiction based on
Lemma~\ref{lem1} and Remark~\ref{u00} leads to the conclusion.


\hfill $\Box$

\bibliographystyle{plain}
\bibliography{biblio}

\end{document}